\documentclass[11pt]{article}
\usepackage{graphicx,wrapfig}
\usepackage{flafter}
\usepackage{leftidx}
\usepackage{mathrsfs}
\usepackage{amssymb,amsmath}
\usepackage{bbding}
\usepackage{fancyhdr}
\usepackage{mathrsfs}
\usepackage{color}
\usepackage{stmaryrd}
\usepackage{amsfonts}
\usepackage{latexsym}
\usepackage{psfrag}
\usepackage{graphicx,subfigure}
\usepackage{fancyhdr,graphicx}
\usepackage{multicol}
\usepackage{dsfont}
\usepackage{bbm}
\usepackage{booktabs}
\usepackage[center]{caption2}
\usepackage{cite}
\usepackage{epstopdf}
\usepackage{booktabs}
\usepackage [latin1]{inputenc}
\usepackage{multirow}
\usepackage{enumerate}
\usepackage{enumitem}
\usepackage{algpseudocode,algorithm,algorithmicx}

\allowdisplaybreaks
\renewcommand{\leq}{\leqslant}

\date{}
\newtheorem{theorem}{Theorem}[section]
\newtheorem{lemma}{Lemma}[section]
\newtheorem{remark}{Remark}
\newtheorem{example}{Example}[section]

\numberwithin{equation}{section}
\newcommand{\zd}{\,\mathrm{d}}

\newcommand{\abs}[1]{\left|#1\right|}

\newcommand{\abst}[1]{|#1|}
\newcommand{\bra}[1]{\left(#1\right)}
\newcommand{\brab}[1]{\big(#1\big)}

\newcommand{\brat}[1]{(#1)}
\newcommand{\kbra}[1]{\left[#1\right]}
\newcommand{\kbrab}[1]{\big[#1\big]}
\newcommand{\kbraB}[1]{\Big[#1\Big]}

\newcommand{\mynorm}[1]{\left\|#1\right\|}

\begin{document}
\title{Adaptive second-order Crank-Nicolson time-stepping schemes for
 time fractional molecular beam epitaxial growth models}
\author{Bingquan Ji\thanks{Department of Mathematics, Nanjing University of Aeronautics and Astronautics,
211101, P. R. China. Bingquan Ji (jibingquanm@163.com).}
\quad Hong-lin Liao\thanks{ORCID 0000-0003-0777-6832; Corresponding author. Department of Mathematics,
Nanjing University of Aeronautics and Astronautics,
Nanjing 211106, P. R. China. Hong-lin Liao (liaohl@csrc.ac.cn)
is supported by a grant 1008-56SYAH18037
from NUAA Scientific Research Starting Fund of Introduced Talent.}
\quad Yuezheng Gong \thanks{Department of Mathematics, Nanjing University of Aeronautics and Astronautics,
Nanjing 210016, P. R. China;
Yuezheng Gong (gongyuezheng@nuaa.edu.cn) is partially supported by the NSFC grant No. 11801269,
 and the NSF grant No. BK20180413 of Jiangsu Province.}
\quad Luming Zhang\thanks{Department of Mathematics, Nanjing University of Aeronautics and Astronautics,
211101, P. R. China. Luming Zhang (zhanglm@nuaa.edu.cn)
is supported by the NSFC grant No. 11571181.}}
\date{}
\maketitle
\normalsize

\begin{abstract}
  Adaptive second-order Crank-Nicolson time-stepping methods using the recent scalar auxiliary variable (SAV) approach
  are developed for the time-fractional Molecular Beam Epitaxial models with Caputo's derivative.
  Based on the piecewise linear interpolation,
  the Caputo's fractional derivative is approximated by a novel second-order formula,
  which is naturally suitable for a general class of nonuniform meshes
  and essentially preserves the positive semi-definite property of integral kernel.
  The resulting Crank-Nicolson SAV time-stepping schemes are unconditional energy stable
  on nonuniform time meshes,
  and are computationally efficient in multiscale time simulations when combined with adaptive time steps,
  such as are appropriate for accurately resolving the intrinsically initial singularity of solution
  and for efficiently capturing fast dynamics away from the initial time.
  Numerical examples are presented to show the effectiveness of our methods.\\
  \noindent{\emph{Keywords}:}\;\; Time-fractional molecular beam epitaxial;  novel L1-type formula;
  scalar auxiliary variable approach; unconditional energy stability; adaptive time-stepping strategy

  \noindent{\bf AMS subject classiffications.}\;\; 35Q99, 65M06, 65M12, 74A50
\end{abstract}
\section{Introduction}
The Molecular Beam Epitaxial (MBE) growth models, in recent years, have become a powerful new technique in material science, such as making compound semi-conductor manufacture with great precision and high purity.
Also, this technique is broadly applied to investigate thin-film deposition of single crystal.
Roughly speaking, the mathematical models used in previous works to study dynamics of the MBE growth process can be classified into three broad categories known as:
atomistic models that are performed using the form of molecular dynamics
\cite{Clarke1987Origin};
continuum models in the form of partial differential systems
\cite{Villain1991Continuum};
and hybrid models
\cite{Gyure1998Level}
which can be regarded as a compromise in the light of the models mentioned above.

In current work, we are interested in the continuum model for the evolution of the MBE growth that is derived by energy variational strategy.
More precisely, given the effective free energy of the model $E[\phi]$
and associated with $L^{2}$ gradient flow, the height evolution model could be written as
\begin{align}
\partial_{t}\phi
=-M\mu,
\end{align}
in which $\phi$ is a scaled height function of a thin film, positive constant $M$ is the mobility coefficient
and $\mu=\frac{\delta{E}}{\delta\phi}$ is the variational derivative of the free energy $E$.
One intrinsic property of above system is the energy dissipation law, that is,
\begin{align}\label{Integer-Energy-Decay-Law}
\frac{\zd{E}}{\zd{t}}
=\bra{\frac{\delta{E}}{\delta\phi},
\frac{\partial\phi}{\partial{t}}}
=-\frac{1}{M}\mynorm{\partial_{t}\phi}^{2}
\leq{0},
\end{align}
where the $L^{2}$ inner product is defined by
$\bra{f,g}=\int_{\Omega}fg\zd{\mathbf{x}}$, and the $L^{2}$ norm $\mynorm{f}=\sqrt{\int_{\Omega}f^2\zd{\mathbf{x}}}$, for all $f,g\in{L}^{2}(\Omega)$.
Notice that the periodic boundary condition or any other proper boundary that can satisfy the flux free condition at the boundary
$\partial_{\mathbf{n}}\phi|_{\partial\Omega}=0$ and
$\partial_{\mathbf{n}}\Delta\phi|_{\partial\Omega}=0$ are chose to ensure that the boundary integrals resulted during the integration by parts vanish, where $\mathbf{n}$ is the outward normal on the boundary.


There are a great amount of works contributed to the investigation of numerical approximations for the solution of the phase field models, for instance the so-called convex splitting technique
\cite{Shen2012Second}
and
the stabilized semi-implicit method
\cite{Xu2006Stability}.
More recently, there are two novel strategies that are used to design second-order, unconditionally energy stable numerical schemes to solve the phase field models: the invariant energy quadratization (IEQ) method \cite{Yang2017numerical} and the scalar auxiliary variable (SAV) approach \cite{Cheng2019Highly}.
The latter is inspired by the former while leads to numerical schemes that only decoupled equations with constant coefficients need to be solved at each time step.
Notwithstanding, the common goal of IEQ and SAV methods is to transform the primitive system into a new equivalent system with a quadratic energy functional and the corresponding modified energy dissipation property.
Rigorous analysis and comparisons between IEQ and SAV might be out of this article's scope, we refer to
\cite{Yang2017numerical,Gong2019Energy,Cheng2019Highly}
for more details.

In comparison with the bright achievement of classical phase field models, there are many researches on building fractional phase field models to better model the anomalously diffusive effects.
For instance,  the time, space and time-space fractional Allen-Cahn equations were suggested
{\cite{Zheng2017A, Hou2017Numerical,Liu2018Time,Zhao2019On}}
to accurately describe anomalous diffusion problems.
Li et al. {\cite{Zheng2017A}} investigated a space-time fractional Allen-Cahn phase field model
that describes the transport of the fluid mixture of two immiscible fluid phases.
They concluded that the alternative model could provide more accurate description of anomalous diffusion processes
and sharper interfaces than the classical model.
Hou et al. {\cite{Hou2017Numerical}} showed that the space-fractional Allen-Cahn equation could be
viewed a $L^{2}$ gradient flow for the fractional analogue version of Ginzburg-Landau free energy function.
Also, the authors proved the energy decay property and the maximum principle of continuous problem.
Tang et al. \cite{Tang2018On} proved the time-fractional phase field models indeed admit an energy dissipation law of an integral type.
Meanwhile, they applied the uniform L1 formula to construct a class of finite difference schemes,
which can inherit the theoretical energy dissipation property.
Along the numerical front, Zhao et al. \cite{Liu2018Time,Zhao2019On} studied a series of the
time-fractional phase field models numerically, covering the time-fractional Cahn-Hilliard equation with different types of variable mobilities and  time-fractional molecular beam epitaxy model.
The considerable numerical evidences
indicate that the effective free energy or roughness of the time-fractional phase field models
during coarsening obeys a similar power scaling law as the integer ones,
where the power is linearly proportional to the fractional index $\alpha$.
In other words, the main difference between the time-fractional phase field models and integer ones lie in the timescales of coarsening.

The multi-scale nature of time-fractional phase field models prompts us to construct reliable time-stepping methods
on general nonuniform meshes. In this paper, nonuniform time-stepping schemes
are investigated for the time-fractional MBE model
\begin{align}\label{Problem-1}
\partial_{t}^{\alpha}\phi
=-M\mu,
\end{align}
where the notation $\partial_{t}^{\alpha}:={}_{0}^{C}\!D_{t}^{\alpha}$ in {\eqref{Problem-1}} denotes
the fractional Caputo derivative of order $\alpha$ with respect to $t$,
\begin{align}\label{CaputoDef}
(\partial_{t}^{\alpha}v)(t)
:=(\mathcal{I}_{t}^{1-\alpha}v')(t)=\int_{0}^{t}\omega_{1-\alpha}(t-s)v'(s)\zd{s},\quad 0<\alpha<1,
\end{align}
involving the fractional  Riemann-Liouville integral $\mathcal{I}_{t}^{\beta}$ of order $\beta>0$, that is,
\begin{align}
(\mathcal{I}_{t}^{\beta}v)(t)
:=\int_{0}^{t}\omega_{\beta}(t-s)v(s)\zd{s},\quad\text{where}\quad  \omega_{\beta}(t):=t^{\beta-1}/\Gamma(\beta).
\end{align}
Specifically, in comparison to the decay property \eqref{Integer-Energy-Decay-Law} of the classical model, the energy dissipation law of the time-fractional MBE model {\eqref{Problem-1}}, see also \cite{Tang2018On}, is given by
\begin{align}\label{Frac-Energy-Decay-Law}
E\kbra{\phi(T)}
-E\kbra{\phi(0)}
=-\frac{1}{M}
\int_{\Omega}\int_{0}^{T}\brab{\partial_{t}^{\alpha}\phi}\partial_{t}\phi\zd{t}\zd{\mathbf{x}}
\leq{0}.
\end{align}
To our knowledge, there are few results in the literature on the discrete energy decay laws
of numerical approaches for the time-fractional phase field models,
especially on nonuniform time meshes.
One of our interests in this paper is to build  nonuniform time-stepping methods preserving the energy dissipation law
of the problem {\eqref{Problem-1}} in discrete sense.

We consider the nonuniform time levels $0=t_{0}<t_{1}<\cdots<t_{k-1}<t_{k}<\cdots<t_{N}=T$ with the time-step sizes $\tau_{k}:=t_{k}-t_{k-1}$ for $1\leq{k}\leq{N}$ and the maximum time-step size $\tau:=\max_{1\leq{k}\leq{N}}\tau_{k}$. Also, let the local time-step ratio $\rho_k:=\tau_k/\tau_{k+1}$ and the maximum step ratio $\rho:=\max_{k\geq 1}\rho_k$. Given a grid function $\{v^{k}\}$, put $\triangledown_{\tau}v^{k}:=v^{k}-v^{k-1}$, $\partial_{\tau}v^{k-\frac12}:=\triangledown_{\tau}v^{k}/\tau_k$ and
$v^{k-\frac{1}{2}}:=(v^{k}+v^{k-1})/2$ for $k\geq{1}$. Always,
let $(\Pi_{1,k}v)(t)$ denote the linear interpolant of a function $v(t)$ at two nodes $t_{k-1}$ and $t_{k}$,
and define a piecewise linear approximation
\begin{align}\label{linear interpolation}
\Pi_{1}v:=\Pi_{1,k}v\quad
\text{so that}\quad
(\Pi_{1}v)'(t)=\partial_{\tau}v^{k-\frac12}\quad
\text{for $t_{k-1}<{t}\leq t_{k}$ and $k\geq1$}.
\end{align}

\begin{figure}[htb!]
\centering
\includegraphics[width=2.0in]{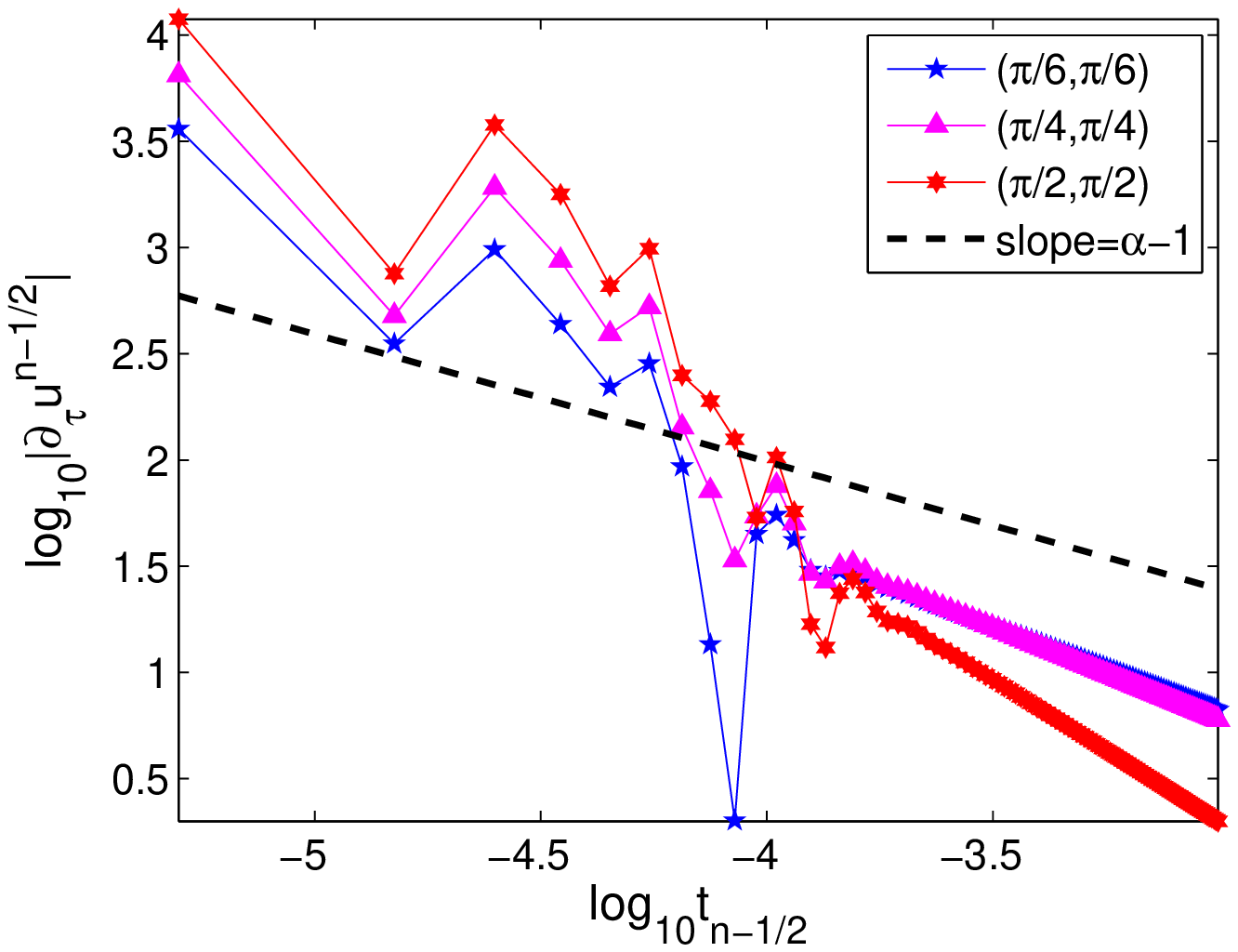}
\includegraphics[width=2.0in]{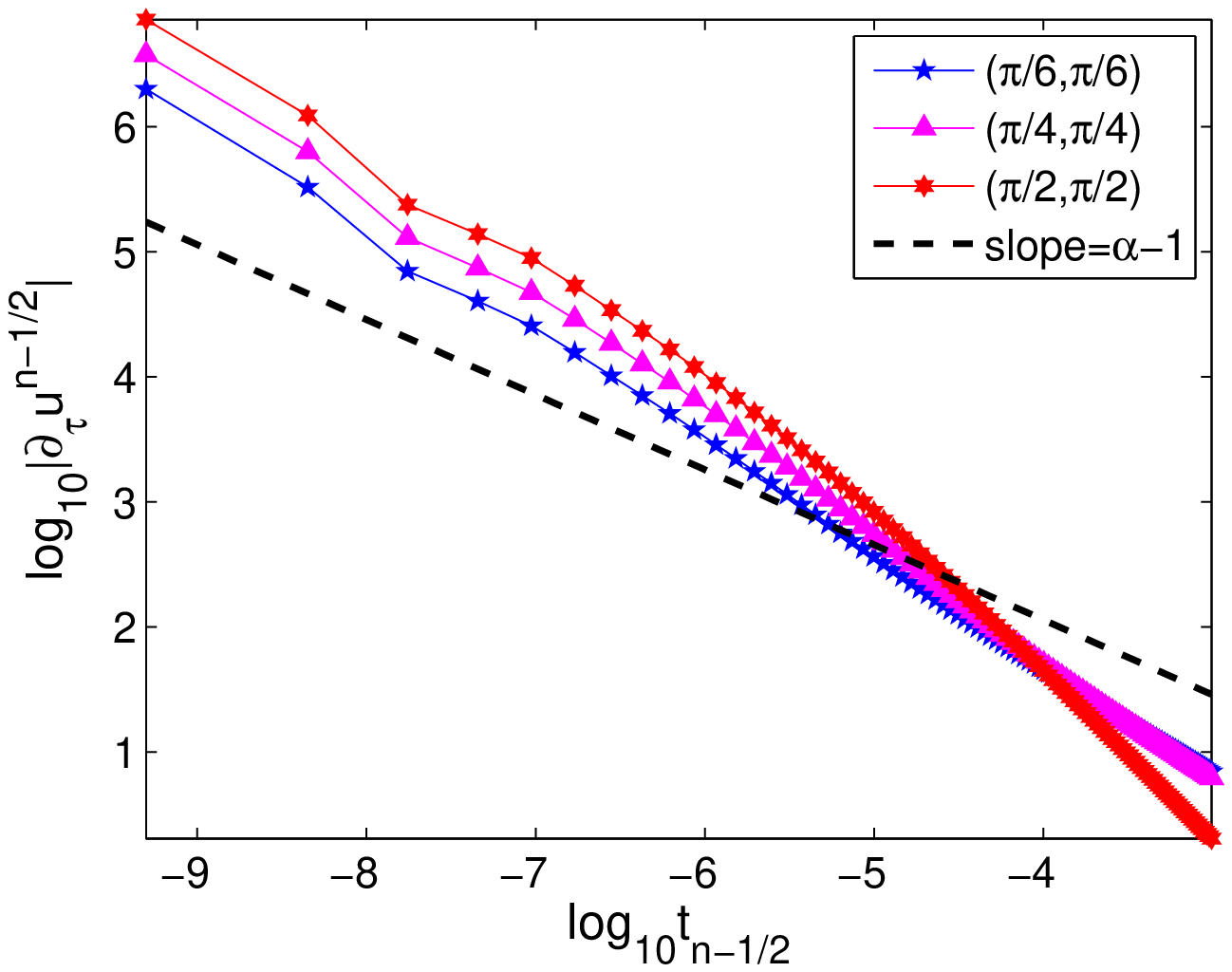}
\includegraphics[width=2.0in]{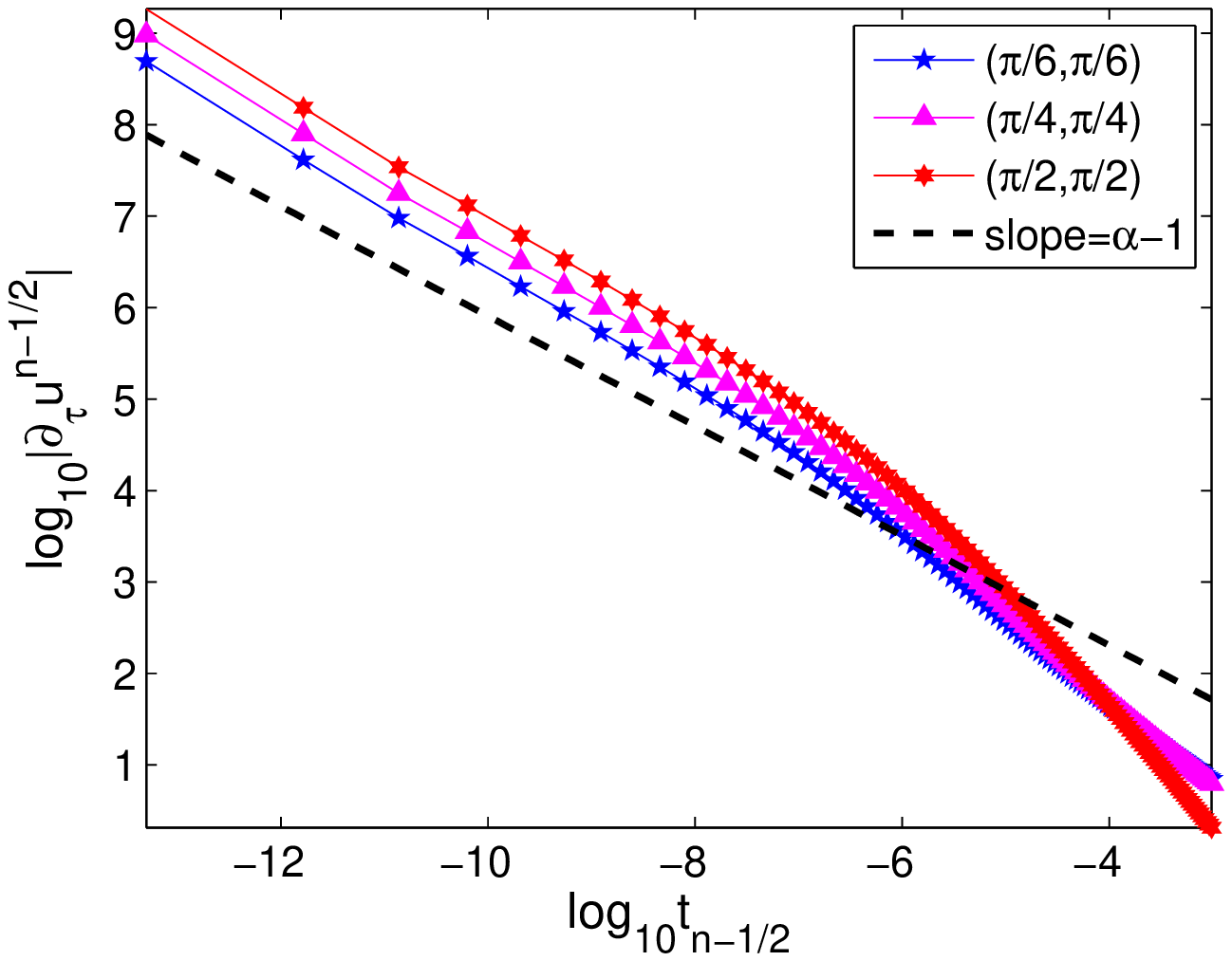}
\caption{The log-log plot of the difference quotient $\partial_{\tau}u^{k-\frac12}$ versus time
for problem {\eqref{Problem-1}} with fractional order $\alpha=0.4$ and $\gamma=1,\,3,\,5$ (from left to right), respectively.}
\label{Initial-Singularity-alpha-04}
\end{figure}

As an essential mathematical feature of linear and nonlinear subdiffusion problems
including the time-fractional MBE model {\eqref{Problem-1}},
the solution always lacks the smoothness near the initial time
although it would be smooth away from  $t=0$, see \cite{Jin2016An,JinLiZhou:2017}.
Actually, in any numerical methods for solving
time-fractional diffusion equations, a key consideration is the initial singularity of solution, see the recent works \cite{Liao2018Sharp,Liao2018Unconditional,Liao2018second}.
More directly, we apply a new L1-type  formula to the time-fractional problem {\eqref{Problem-1}}, see more details in subsection 3.1 and Example \ref{Time-Adaptive-Test-CN-SAV}.
Figure \ref{Initial-Singularity-alpha-04}  plots the discrete time derivative $\partial_{\tau}u^{k-\frac12}$
near $t=0$ on the graded mesh $t_k=(k/N)^{\gamma}$ when $\alpha=0.4$.
They suggest that
\[
\log|u_{t}(\mathbf{x},t)|
\approx(\alpha-1)\log(t)+C(\mathbf{x})\quad\text{as $t\rightarrow0$.}
\]
It says that the solution possesses weakly singularity like $u_{t}=O(t^{\alpha-1})$ near the initial time, which can be alleviated by using the graded meshes.
Thus, the second interest of this paper is to resolve the essentially weak singularity
in the equation {\eqref{Problem-1}} by refining time mesh. Actually, we will show that
the graded mesh can recover the optimal time accuracy when the solution $u$ is non-smooth near $t=0$.

In the next section, we reformulate the time-fractional MBE model {\eqref{Problem-1}} by using SAV technique, which provides an elegant platform for constructing energy stable schemes.
Although the L1 formula was applied in \cite{Tang2018On} to handle the time-fractional phase filed models
and the discrete energy decay laws were also established on uniform time mesh,
we are not able to exploit any discrete energy dissipation laws under the adaptive time-steps
and leave it as an open problem (see Remark \ref{L1-Quadratic-Form}).
In section 3, this open problem guides us to develop a novel numerical Caputo derivative \eqref{New-L1-Formula}, called L1$^{+}$ formula,
which preserves inherently the positive semi-definite property of integral kernel, see \eqref{New-L1-Positive}.
This new formula is as simple as the classical L1 formula because it is also constructed by the piecewise linear interpolation;
however it can achieve second-order accurate for all fractional orders $\alpha\in(0,1)$
and substantially improves the accuracy of L1 formula, especially when the fractional order $\alpha\rightarrow1$.
Applying the new formula to the time-fractional MBE model \eqref{Problem-1},
we obtain two Crank-Nicolson SAV schemes \eqref{Slope-CN-SAV-1}-\eqref{Slope-CN-SAV-2} and
\eqref{No-CN-SAV-1}-\eqref{No-CN-SAV-2}, respectively.
Thus it is straightforward to confirm that both of them preserve the discrete energy dissipation law well,
see Theorems \ref{Slope-CN-SAV-Decay-Law}-\ref{No-CN-SAV-Decay-Law},
so that they are unconditionally stable in the energy norm. At last,
we also apply the sum-of-exponentials (SOE) technique to develop a fast L1$^{+}$ algorithm in section 4.
By using an accuracy criterion based adaptive time-stepping strategy, extensive numerical experiments are curried out
to show the effectiveness of our numerical approaches and to support our analysis.

In summary, the main contributions of this paper for the time-fractional MBE model \eqref{Problem-1} are the following: suggest a simple L1$^{+}$ formula of Caputo derivative having second-order accuracy for any fractional order $\alpha\in(0,1)$;
apply it to build two Crank-Nicolson SAV schemes preserving the discrete energy dissipation law; and develop a fast L1$^{+}$ algorithm
incorporating adaptive time steps to speed up the long-time, multi-scale simulations.

\section{Equivalent PDE systems}
Let $\varepsilon>0$ is a positive constant and $F\bra{\mathbf{v}}$ is a nonlinear, smooth function of its argument $\mathbf{v}\in\mathbb{R}^{d}\,(d=1,2,3)$.
The Ehrlich-Schwoebel energy $E[\phi]$ is given by
\begin{align}
E[\phi]=\int_{\Omega}\bra{\frac{\varepsilon^{2}}{2}\abs{\Delta\phi}^{2}
+ F\bra{\nabla\phi}}\zd{\mathbf{x}}\,.
\end{align}
There are two popular choices of the nonlinear bulk potential:
the double well potential, $F\bra{\nabla\phi}=\frac{1}{4}\brat{\abst{\nabla\phi}^{2}-1}^{2}$,
for the case of slope selection model;
and the logarithmic potential, $F\bra{\nabla\phi}=-\frac{1}{2}\ln\brat{1+\abst{\nabla\phi}^{2}}$,
for the case without slope selection model.
Correspondingly, the governing equation of {\eqref{Problem-1}} for the height function $\phi$ reads:
\begin{align}\label{Problem-2}
\partial_{t}^{\alpha}\phi
=-M\frac{\delta{E}\bra{\phi}}{\delta\phi}
=-M
\bra{
\varepsilon^{2}\Delta^{2}\phi
+f\bra{\nabla\phi}
},
\end{align}
in which the case of slope selection model
\begin{align}\label{Slope-Model}
f\bra{\nabla\phi}
=-\nabla\cdot
\brab{
\brab{\abst{\nabla\phi}^{2}-1}\nabla\phi
},
\end{align}
while the case without slope selection model
\begin{align}\label{No-Slope-Model}
f\bra{\nabla\phi}
=\nabla\cdot
\brab{
\frac{1}{1+\abst{\nabla\phi}^{2}}\nabla\phi
}.
\end{align}
In what follows, we denote the time fractional MBE model \eqref{Problem-2} with and without slope selection by ``the Slope-Model" and ``the No-Slope-Model", respectively, for simplicity.
Additionally, boundary conditions are set to be periodic so as not to complicate the analysis with unwanted details.
We then reformulate the MBE models into the equivalent PDE systems by borrowing the ideas from the IEQ and SAV methods \cite{Yang2017numerical,Gong2019Energy,Cheng2019Highly}.

\subsection{The Slope-Model and its equivalent system}
We introduce a scalar auxiliary function $u(t)$ in term of original variable $\phi$ given by
\begin{align}\label{Slope-New-Variable}
u(t)
=\sqrt{\int_{\Omega}
\frac{1}{4}\bra{\abs{\nabla\phi}^{2}-1-\beta}^{2}\zd{\mathbf{x}}+C_{0}},
\end{align}
where $C_{0}>0$ is constant that ensures the radicand positive.
Compared with those in \cite{Yang2017numerical,Cheng2019Highly}, the scalar auxiliary variable $u(t)$ adds
an artificial parameter $\beta$ to regularize the numerical approach.
As a consequence, the free energy of the Slope-Model is transformed into a quadratic form
\begin{align}\label{Slope-Equivalent-Energy}
E\kbra{\phi,u}
=\int_{\Omega}\bra{
\frac{\varepsilon^{2}}{2}\abs{\Delta\phi}^{2}
+\frac{\beta}{2}\abs{\nabla\phi}^{2}}
\zd{\mathbf{x}}
+u^{2}
-\brab{\frac{\beta}{2}+\frac{\beta^{2}}{4}}\abs{\Omega}-C_{0}.
\end{align}
Correspondingly, the Slope-Model can be reformulated to the following equivalent form
\begin{align}
\partial_{t}^{\alpha}\phi
&=-M
\bra{
\varepsilon^{2}\Delta^{2}\phi
-\beta\Delta\phi
-U\bra{\phi}u},\label{Slope-SAV-1}\\
u_{t}
&=-\frac{1}{2}\int_{\Omega}U\bra{\phi}\partial_{t}\phi\zd{\mathbf{x}},\label{Slope-SAV-2}
\end{align}
where the expression of notation $U\bra{\phi}$ given by
\begin{align}
U\bra{\phi}
=\frac{\nabla\cdot
\brab{
\brab{\abs{\nabla\phi}^{2}-1-\beta}\nabla\phi
}}
{\sqrt{\int_{\Omega}
\frac{1}{4}\bra{\abs{\nabla\phi}^{2}-1-\beta}^{2}\zd{\mathbf{x}}+C_{0}}}.
\end{align}
We should note in passing that
the new system is subjected to the initial conditions
\begin{align}
\phi\bra{\mathbf{x},0}
=\phi_{0}\bra{\mathbf{x}}\quad
\text{and}\quad
u\bra{0}
=u\bra{\phi_{0}\bra{\mathbf{x}}},
\end{align}
and the boundary conditions are same as the primitive model.
By taking the $L^{2}$ inner product of \eqref{Slope-SAV-1} with $\phi$, of \eqref{Slope-SAV-2} with $u$, and then making time integrations on both sides, we see the equivalent system preserves the energy dissipation law
\begin{align}
E\kbra{\phi(T),u(T)}
-E\kbra{\phi(0),u(0)}
=-\frac{1}{M}
\int_{\Omega}\int_{0}^{T}\brab{\partial_{t}^{\alpha}\phi}\partial_{t}\phi\zd{t}\zd{\mathbf{x}}
\leq{0}.\nonumber
\end{align}
The non-positive of the right part of above equality is determined by \cite[Lemma 2.1]{Tang2018On}.

\subsection{The No-Slope-Model and its equivalent system}
For the No-Slope-Model, we introduce a scalar auxiliary function $v(t)$ in term of original variable $\phi$ as follows
\begin{align}
v(t)
=\sqrt{\int_{\Omega}
\bra{
\frac{1}{2}\ln\brab{1+\abs{\nabla\phi}^{2}}
+\frac{\beta}{2}\abs{\nabla\phi}^{2}}\zd{\mathbf{x}}+C_{0}},
\end{align}
where $\beta$ and $C_{0}$ are similar to the previous ones. Therefore, the free energy of the No-Slope-Model could be rewritten into
\begin{align}
E\kbra{\phi,v}=\int_{\Omega}
\bra{\frac{\varepsilon^{2}}{2}\abs{\Delta\phi}^{2}
+\frac{\beta}{2}\abs{\nabla\phi}^{2}}
\zd{\mathbf{x}}
-v^{2}+C_{0}.
\end{align}
We then could rewrite the No-Slope-Model as an equivalent form
\begin{align}
\partial_{t}^\alpha\phi
&=-M
\bra{
\varepsilon^{2}\Delta^{2}\phi
-\beta\Delta\phi
+V(\phi)v},\label{No-Slope-SAV-1}\\
v_{t}
&=\frac{1}{2}\int_{\Omega}
V\bra{\phi}\partial_{t}\phi\zd{\mathbf{x}},\label{No-Slope-SAV-2}
\end{align}
in which
\begin{align}
V(\phi)=
\frac{\nabla\cdot\bra{\brab{\frac{1}{1+\abs{\nabla\phi}^{2}}+\beta}\nabla\phi}}
{\sqrt{\int_{\Omega}
\bra{
\frac{1}{2}\ln\brab{1+\abs{\nabla\phi}^{2}}
+\frac{\beta}{2}\abs{\nabla\phi}^{2}}\zd{\mathbf{x}}+C_{0}}},
\end{align}
with the following initial conditions
\begin{align}
\phi\bra{\mathbf{x},0}
=\phi_{0}\bra{\mathbf{x}}\quad
\text{and}\quad
v(0)
=v\bra{\phi_{0}\bra{\mathbf{x}}}.
\end{align}
Similarly, the new system admits the following energy dissipation law
\begin{align}
E\kbra{\phi(T),v(T)}
-E\kbra{\phi(0),v(0)}
=-\frac{1}{M}
\int_{\Omega}\int_{0}^{T}\brab{\partial_{t}^{\alpha}\phi}\partial_{t}\phi\zd{t}\zd{\mathbf{x}}
\leq{0}.\nonumber
\end{align}

\section{Novel L1$^{+}$ formula and Crank-Nicolson SAV schemes}
The well-known L1 formula of Caputo derivative is given by
\begin{align}\label{L1-Formula}
(\partial_{\tau}^{\alpha}v)^{n}
:=\int_{t_{0}}^{t_{n}}
\omega_{1-\alpha}(t_{n}-s)(\Pi_{1}v)'(s)\zd{s}
=\sum_{k=1}^{n}a_{n-k}^{(n)}\triangledown_{\tau}v^{k},
\end{align}
where the corresponding discrete convolution kernels $a_{n-k}^{(n)}$ are defined by
\begin{align}\label{L1-Formula-Coefficient}
a_{n-k}^{(n)}
:=\frac{1}{\tau_{k}}\int_{t_{k-1}}^{t_{k}}\omega_{1-\alpha}(t_n-s)\zd{s}\quad\text{for $1\leq{k}\leq{n}.$}
\end{align}
Obviously, the discrete L1 kernels $a_{n-k}^{(n)}$ are positive and decreasing, see also \cite{Liao2018Sharp,Liao2018Unconditional},
\begin{align}\label{L1-Coefficient-Estimate}
a_{n-k}^{(n)}>0 \quad\text{and}\quad a_{n-k-1}^{(n)}>a_{n-k}^{(n)}
\quad\text{for $1\leq{k}\leq{n-1}.$}
\end{align}

Based on the above L1 formula on uniform time mesh, a linearized scheme by using the stabilized technique
via a stabilized term $S(\Delta\phi^{n}-\Delta\phi^{n-1})$ for a properly large scalar parameter $S>0$ is presented. We refer to \cite{Tang2018On} for more details.
The resulting stabilized semi-implicit scheme for the problem \eqref{Problem-2} reads
\begin{align}\label{Stabilized-Scheme}
\bra{\partial_{\tau}^\alpha\phi}^{n}
=-M\bra{
\varepsilon^{2}\Delta^{2}\phi^{n}
+f\bra{\nabla\phi^{n-1}}
}
+S(\Delta\phi^{n}-\Delta\phi^{n-1}).
\end{align}
It is to mention that, when the time mesh is uniform such that the discrete L1 kernels
\begin{align*}
a_{n-k}^{(n)}=a_{n-k}
=\frac{1}{\tau^{\alpha}}\kbra{\omega_{2-\alpha}(n-k+1)-\omega_{2-\alpha}(n-k)}\quad\text{for $1\leq{k}\leq{n},$}
\end{align*}
the semi-implicit scheme \eqref{Stabilized-Scheme}
inherits a discrete energy dissipation law, see more details in \cite[Theorem 3.3]{Tang2018On}.
As seen, the proof of discrete version of energy dissipation law
\eqref{Frac-Energy-Decay-Law} relies on
the positive semi-definite property of a discrete quadratic form $\sum_{k=1}^nw_k\sum_{j=1}^ka_{k-j}w_j\ge0$.

\begin{remark}\label{L1-Quadratic-Form}
It seems rather difficult to extend the positive semi-definite property to a general class of nonuniform meshes.
More precisely, we are not able to verify the positive semi-definite property of the following quadratic form
(by taking $w_k=\triangledown_{\tau}v^k$)
\begin{align}\label{Positive-Quadratic-Form}
\sum_{k=1}^n\triangledown_{\tau}v^{k}(\partial_{\tau}^{\alpha}v)^{k}
=\sum_{k=1}^nw_k\sum_{j=1}^ka_{k-j}^{(k)}w_j\ge0.
\end{align}
\end{remark}

The open problem in Remark \ref{L1-Quadratic-Form} motivates us to design a novel discrete Caputo formula such that
it naturally possesses the energy decay law \eqref{Frac-Energy-Decay-Law} in the discrete sense
for the time-fractional MBE model.
Alternatively, the novel discrete Caputo formula should inherit the positive semi-definite property
of a quadratic form like \eqref{Positive-Quadratic-Form}.

Fortunately, as pointed out early in \cite{Mclean1996Discretization,William2007A},
we know that the
weakly singular kernel $\omega_{1-\alpha}$ is positive semi-definite, that is,
\begin{align}\label{Continu-Positive-Kernel}
\mathcal{I}_{t}^{1}(w\,\mathcal{I}_{t}^{1-\alpha}w)(t)
&=\int_{0}^{t}w(\mu)\zd{\mu}\int_{0}^{\mu}\omega_{1-\alpha}(\mu-s)w(s)\zd{s}\nonumber\\
&=\frac{1}{2}\int_{0}^{t}\int_{0}^{t}\omega_{1-\alpha}(|\mu-s|)w(\mu)w(s)\zd{\mu}\zd{s}\geq{0}
\end{align}
for $t>0$ and any $w\in{C}[0,T]$. Actually, one has $(-1)^k\omega_{1-\alpha}^{(k)}(t)\geq{0}$ for any integer $k\geq1$
and $\int_{0}^{\infty}\omega_{1-\alpha}(t)\cos(\theta{t})\zd{t}\geq{0}$ for any $\theta>0$, then the Plancherel's theorem implies
the positive semi-definite property \eqref{Continu-Positive-Kernel}.
We will see that a discrete counterpart of \eqref{Continu-Positive-Kernel} yields a discrete Caputo approximation preserving
the desired energy dissipation property \eqref{Frac-Energy-Decay-Law} when it is applied to
time-fractional MBE model \eqref{Problem-2}.

\subsection{The L1$^{+}$ formula}
The L1$^{+}$ formula for the Caputo derivative \eqref{CaputoDef} is defined at time $t=t_{n-\frac{1}{2}}$ as follows
\begin{align}\label{New-L1-Formula}
(\partial_{\tau}^{\alpha}v)^{n-\frac{1}{2}}
:=\frac{1}{\tau_{n}}\int_{t_{n-1}}^{t_{n}}
\int_{0}^{t}\omega_{1-\alpha}(t-s)(\Pi_{1}v)'(s)\zd{s}\zd{t}
=\sum_{k=1}^{n}\bar{a}_{n-k}^{(n)}\triangledown_{\tau}v^{k}\quad \text{for $n\ge1$,}
\end{align}
where the discrete convolution kernels $\bar{a}_{n-k}^{(n)}$ are defined by
\begin{align}\label{New-L1-Coeff}
\bar{a}_{n-k}^{(n)}
:=\frac{1}{\tau_{n}\tau_{k}}\int_{t_{n-1}}^{t_{n}}
\int_{t_{k-1}}^{\min\{t,t_{k}\}}\omega_{1-\alpha}(t-s)\zd{s}\zd{t}\quad \text{for $1\leq{k}\leq{n}$.}
\end{align}
Obviously, the naturally nonuniform L1$^{+}$ approximation $(\partial_{\tau}^{\alpha}v)^{n-\frac{1}{2}}\approx\frac{1}{\tau_{n}}\int_{t_{n-1}}^{t_{n}}\bra{\partial_{t}^{\alpha}v}(t)\zd{t}$ 
ensures the positive semi-definite property \eqref{Continu-Positive-Kernel}
by taking $w=\Pi_{1}v$, that is,
\begin{align}
\sum_{k=1}^{n}\triangledown_{\tau}v^{k}(\partial_{\tau}^{\alpha}v)^{k-\frac{1}{2}}
&=\sum_{k=1}^{n}\tau_{k}(\Pi_{1,k}v)'(\partial_{\tau}^{\alpha}v)^{k-\frac{1}{2}}\nonumber\\
&=\int_{t_{0}}^{t_{n}}(\Pi_{1}v)'(t)
\int_{0}^{t}\omega_{1-\alpha}(t-s)(\Pi_{1}v)'(s)\zd{s}\zd{t}\nonumber\\
&=\mathcal{I}_{t}^{1}\kbraB{(\Pi_{1}v)'
\mathcal{I}_{t}^{1-\alpha}(\Pi_{1}v)'}(t_{n})\geq{0}\quad \text{for $n\ge1$.}\label{New-L1-Positive}
\end{align}
The definition \eqref{New-L1-Coeff} and the arbitrariness of function $v$ yield the following result.
\begin{lemma}\label{lem: CN quadraticForm}
The discrete convolution kernels $\bar{a}_{n-k}^{(n)}$
in \eqref{New-L1-Coeff} are positive, and for any real sequence $\{w_k\}_{k=1}^n$ with $n$ entries, it holds that
\begin{align*}
\sum_{k=1}^nw_k\sum_{j=1}^k\bar{a}_{k-j}^{(k)}w_j\geq{0}\quad \text{for $n\ge1$.}
\end{align*}
\end{lemma}

\begin{table}[htbp]
\begin{center}
\caption{Numerical accuracy of L1$^{+}$ formula \eqref{New-L1-Formula} with $\sigma=2.5$}\label{New-L1-Error-1}  \vspace*{0.5pt}
\def\temptablewidth{1.0\textwidth}
{\rule{\temptablewidth}{0.5pt}}
\begin{tabular*}{\temptablewidth}{@{\extracolsep{\fill}}ccccccc}
\multirow{2}{*}{$N$} &\multicolumn{2}{c}{$\alpha=0.1$} &\multicolumn{2}{c}{$\alpha=0.5$} &\multicolumn{2}{c}{$\alpha=0.9$}\\
          \cline{2-3}       \cline{4-5}       \cline{6-7}
          &$e(N)$  &Order     &$e(N)$ &Order     &$e(N)$  &Order \\
\midrule
  64      &3.44e-05  &--        &3.39e-05   &--      &2.25e-05    &--\\
  128     &8.61e-06  &2.00      &8.51e-06   &1.99    &5.83e-06    &1.95\\
  256     &2.15e-06  &2.00      &2.13e-06   &1.99    &1.51e-06    &1.95\\
  512     &5.38e-07  &2.00      &5.35e-07   &2.00    &3.88e-07    &1.96\\
\end{tabular*}
{\rule{\temptablewidth}{0.5pt}}
\end{center}
\end{table}
Before applying the L1$^{+}$ formula \eqref{New-L1-Formula} to the time-fractional MBE model {\eqref{Problem-2}},
we show that it is a second-order approximation for the Caputo derivative \eqref{CaputoDef} numerically.
Consider a simple fractional ODE problem
$\partial_{t}^{\alpha}u(t)=f(t)$ for $0<t<1$, we run a Crank-Nicolson-type scheme
$(\partial_{\tau}^{\alpha}u)^{n-\frac{1}{2}}=f(t_{n-\frac12})$
with uniform time-steps $\tau_k=\tau$,
by choosing a smooth solution $u=\omega_{1+\sigma}(t)$ with the regularity parameter $\sigma=2.5$.
The discrete $l^{\infty}$ norm errors $e(N)=\max_{1\leq{n}\leq{N}}|u(t_{n})-u^{n}|$ are listed in Table \ref{New-L1-Error-1}.
It seems that the numerical accuracy of \eqref{New-L1-Formula} is second-order accurate for any fractional order $\alpha\in(0,1)$.
If the solution has an initial singularity, the L1$^{+}$ formula \eqref{New-L1-Formula} can also achieve the second-order accuracy
by properly refining the mesh near $t=0$, see more tests in Example \ref{Accuracy-Test-New-L1}.

\subsection{Crank-Nicolson SAV schemes preserving energy dissipation}
In what follows, we concern only with  the time discretization of the equivalent systems, while the spatial approximation can be diverse, examples as finite difference, finite element or spectral methods.
Integrating the equations \eqref{Slope-SAV-1}-\eqref{Slope-SAV-2} from $t=t_{n-1}$ to $t_{n}$,
respectively, leads to
\begin{align}
\frac{1}{\tau_{n}}\int_{t_{n-1}}^{t_{n}}\partial_{t}^{\alpha}\phi\zd{t}
&=-\frac{M}{\tau_{n}}\int_{t_{n-1}}^{t_{n}}
\bra{
\varepsilon^{2}\Delta^{2}\phi
-\beta\Delta\phi
-U\bra{\phi}u}\zd{t},\\
\frac{1}{\tau_{n}}\int_{t_{n-1}}^{t_{n}}u_{t}\zd{t}
&=-\frac{1}{2\tau_{n}}\int_{t_{n-1}}^{t_{n}}
\int_{\Omega}U\bra{\phi}\partial_{t}\phi\zd{\mathbf{x}}\zd{t}.
\end{align}
Applying the L1$^{+}$ formula \eqref{New-L1-Formula}, the trapezoidal formula, we have the following Crank-Nicolson SAV (CN-SAV) time-stepping  scheme for the Slope-Model
\begin{align}
\bra{\partial_{\tau}^{\alpha}\phi}^{n-\frac{1}{2}}
&=-M
\bra{
\varepsilon^{2}\Delta^{2}\phi^{n-\frac{1}{2}}
-\beta\Delta\phi^{n-\frac{1}{2}}
-U(\hat{\phi}^{n-\frac{1}{2}})u^{n-\frac{1}{2}}
},\label{Slope-CN-SAV-1}\\
\partial_{\tau}u^{n-\frac{1}{2}}
&=-\frac{1}{2}\int_{\Omega}
U(\hat{\phi}^{n-\frac{1}{2}})
\partial_{\tau}\phi^{n-\frac{1}{2}}\zd{\mathbf{x}},\label{Slope-CN-SAV-2}
\end{align}
where $\hat{\phi}^{n-\frac{1}{2}}:=\phi^{n-1}+\triangledown_{\tau}\phi^{n-1}/(2\rho_{n-1})$ is the local extrapolation.

Note that, the construction of L1$^{+}$ formula \eqref{New-L1-Formula}
implies that the above CN-SAV scheme \eqref{Slope-CN-SAV-1}-\eqref{Slope-CN-SAV-2} is naturally suitable for a general class of nonuniform meshes. Moreover,
next result shows that it is unconditionally energy stable.
\begin{theorem}\label{Slope-CN-SAV-Decay-Law}
The CN-SAV scheme \eqref{Slope-CN-SAV-1}-\eqref{Slope-CN-SAV-2} preserves the energy dissipation law,
\begin{align}
E\kbra{\phi^{n},u^{n}}
-E\kbra{\phi^{0},u^{0}}
\leq{0},\quad \text{for}\quad 1\leq{n}\leq{N},
\end{align}
such that it is unconditionally stable, where
\begin{align*}
E\kbra{\phi^{n},u^{n}}
=\int_{\Omega}
\bra{\frac{\varepsilon^{2}}{2}\abs{\Delta\phi^{n}}^{2}
+\frac{\beta}{2}\abs{\nabla\phi^{n}}^{2}}\zd{\mathbf{x}}
+\bra{u^{n}}^{2}-C_{0}.
\end{align*}
\end{theorem}
\begin{proof}
Taking the inner product of \eqref{Slope-CN-SAV-1}-\eqref{Slope-CN-SAV-2} with $\triangledown_{\tau}\phi^{n}$ and $2\tau_{n}u^{n-\frac{1}{2}}$,
respectively, and adding the resulting two identities, we obtain
\begin{align}
-\frac{1}{M}
\bra{\bra{\partial_{\tau}^{\alpha}\phi}^{n-\frac{1}{2}},\triangledown_{\tau}\phi^{n}}
&=
\bra{
\varepsilon^{2}\Delta^{2}\phi^{n-\frac{1}{2}}
-\beta\Delta\phi^{n-\frac{1}{2}},\triangledown_{\tau}\phi^{n}}
+\bra{u^{n}}^{2}-\bra{u^{n-1}}^{2},
\end{align}
which implies that
\begin{align}
E\kbrab{\phi^{k},u^{k}}
-E\kbrab{\phi^{k-1},u^{k-1}}
=-\frac{1}{M}
\bra{(\partial_{\tau}^{\alpha}\phi)^{k-\frac{1}{2}},\triangledown_{\tau}\phi^{k}}\quad \text{for $1\leq{k}\leq{n}$.}
\end{align}
By summing the superscript $k$ from $1$ to $n$, we apply the property \eqref{New-L1-Positive} to derive that
\begin{align*}
E\kbrab{\phi^{n},u^{n}}
-E\kbrab{\phi^{0},u^{0}}
=-\frac{1}{M}\int_{\Omega}\mathcal{I}_{t}^{1}\kbraB{(\Pi_{1}\phi)^{\prime}
\mathcal{I}_{t}^{1-\alpha}(\Pi_{1}\phi)^{\prime}}(t_{n})\zd{\mathbf{x}}\leq{0}\quad \text{for $1\leq{n}\leq{N}.$}
\end{align*}
It completes the proof.
\end{proof}

Analogously, for the No-Slope-Model \eqref{No-Slope-SAV-1}-\eqref{No-Slope-SAV-2}, we present the following CN-SAV time-stepping  scheme
\begin{align}
\bra{\partial_{\tau}^{\alpha}\phi}^{n-\frac{1}{2}}
&=-M
\bra{
\varepsilon^{2}\Delta^{2}\phi^{n-\frac{1}{2}}
-\beta\Delta\phi^{n-\frac{1}{2}}
+V(\hat{\phi}^{n-\frac{1}{2}})v^{n-\frac{1}{2}}},\label{No-CN-SAV-1}\\
\partial_{\tau}v^{n-\frac{1}{2}}
&=\frac{1}{2}\int_{\Omega}V(\hat{\phi}^{n-\frac{1}{2}})
\partial_{\tau}\phi^{n-\frac{1}{2}}\zd{\mathbf{x}}.\label{No-CN-SAV-2}
\end{align}
For the above numerical scheme, we shall only state the energy dissipation law below, as their proofs are essentially the same as Theorem \ref{Slope-CN-SAV-Decay-Law}.
\begin{theorem}\label{No-CN-SAV-Decay-Law}
The CN-SAV scheme \eqref{No-CN-SAV-1}-\eqref{No-CN-SAV-2} satisfies the following discrete energy dissipation law,
\begin{align}
E\kbra{\phi^{n},v^{n}}
-E\kbra{\phi^{0},v^{0}}
\leq{0},
\end{align}
where
\begin{align}
E\kbra{\phi^{n},v^{n}}
=\int_{\Omega}\bra{\frac{\varepsilon^{2}}{2}\abs{\Delta\phi^{n}}^{2}
+\frac{\beta}{2}\abs{\nabla\phi^{n}}^{2}}\zd{\mathbf{x}}
-\bra{v^{n}}^{2}+C_{0}.
\end{align}
\end{theorem}
\subsection{Further notes on L1$^{+}$ formula}

\begin{remark}[Multi-term and distributed time-fractional problems]
\label{Distributed-order-etc}
As seen in Table  \ref{New-L1-Error-1} and more tests for Example
presented in next section,
the accuracy of the L1$^{+}$ formula \eqref{New-L1-Formula} is dependent on the regularity of solution,
but would be independent of the fractional order $\alpha\in(0,1)$.
It is quite different from some exiting numerical approaches, such as L1, Alikhanov \cite{Alikhanov2015A,Liao2016A}
and BDF2-like formulas \cite{Gao2014A,LvXu2016,Liao2016Stability}, of which the consistency errors are dependent on the fractional order $\alpha$.
This feature is attractive for further applications in developing second-order approximations
for multi-term and distributed-order fractional diffusion equations. For an example,
consider a simple multi-term fractional diffusion problem,
\[
\sum_{i=1}^mw_{i}\,\partial_t^{\alpha_{i}}u=\Delta u+f,
\]
where $0<\alpha_{i}<1$ for $1\leq i\leq m$, and  $w_{i}$ is the corresponding weights.
One can construct the following second-order Crank-Nicolson-type time-stepping scheme
\[
\sum_{i=1}^mw_{i}(\partial_{\tau}^{\alpha_{i}}v)^{n-\frac{1}{2}}=\Delta u^{n-\frac{1}{2}}+f^{n-\frac{1}{2}}\quad\text{for $1\leq n\leq N$}.
\]
Obviously, the L1$^{+}$ formula \eqref{New-L1-Formula} would be also useful in
approximating the distributed-order Caputo
derivative  since it can be approximated by certain multi-term
derivative via some proper quadrature rule, see \cite{Liao2016Stability}.
\end{remark}

The suggested L1$^{+}$ formula \eqref{New-L1-Formula} seems very promising in further applications
for other time-fractional field phase models and other time-fractional differential equations.
They make the rigorous theoretical analysis on consistency, stability and convergence very important,
especially on a general class of nonuniform meshes.

However,
the established theory \cite{Liao2018Sharp,Liao2018discrete,Liao2018Unconditional,Liao2018second}
for nonuniform L1 and L1-2$_{\sigma}$ (Alikhanov)  formulas can not be applied to the L1$^{+}$ formula directly
because the corresponding discrete convolution kernels $\bar{a}_{n-k}^{(n)}$ in \eqref{New-L1-Coeff} do not have the uniform monotonicity like \eqref{L1-Formula-Coefficient}.
Actually, the definition \eqref{New-L1-Coeff} and the integral mean-value theorem
yield the following result.
\begin{lemma}\label{New-L1-Coeff-Relation}
The positive discrete kernels $\bar{a}_{n-k}^{(n)}$ in \eqref{New-L1-Coeff} fulfill
\[
\bar{a}_{0}^{(n)}=\frac{1}{\Gamma(3-\alpha)\tau_{n}^{\alpha}},\quad
\bar{a}_{1}^{(n)}>\bar{a}_{2}^{(n)}>\cdots>\bar{a}_{n-1}^{(n)}>0\quad \text{for $n\ge2$}.
\]
\end{lemma}

Notwithstanding, it ought to be emphasized that
\begin{align*}
\bar{a}_{0}^{(n)}-\bar{a}_{1}^{(n)}
=\frac{1}{\Gamma(3-\alpha)\tau_{n}^{\alpha}\rho_{n-1}}
\brab{1+\rho_{n-1}+\rho_{n-1}^{2-\alpha}-(1+\rho_{n-1})^{2-\alpha}}.
\end{align*}
It is easily seen that $\bar{a}_{0}^{(n)}<\bar{a}_{1}^{(n)}$ as $\alpha\rightarrow{0}$ and $\bar{a}_{0}^{(n)}>\bar{a}_{1}^{(n)}$ as $\alpha\rightarrow{1}$. So the value of $\bar{a}_{0}^{(n)}-\bar{a}_{1}^{(n)}$ may change the sign when the fractional order $\alpha$ varies over $(0,1)$. At the same time, this situation is no worse than the case in the BDF2-like formulas \cite{LvXu2016,Liao2016Stability}
in which the second kernel would be negative when $\alpha\rightarrow{1}$, see more details
and a potential remedy technique in \cite[Remark 6]{Liao2018discrete}.
The theoretical investigations,
including the consistency, stability and convergence, of nonuniform L1$^{+}$ formula \eqref{New-L1-Formula}
 will be addressed in a separate technical report.


\section{Numerical algorithms and examples}
\subsection{A fast version of L1$^{+}$ formula}
It is evident that the approximations \eqref{L1-Formula} or \eqref{New-L1-Formula} are prohibitively expensive for long time simulations due to the long-time memory. Therefore, to reduce the computational cost and storage requirements,
we apply the sum-of-exponentials (SOE) technique to speed up the evaluation of the L1$^{+}$ formula.
A core result is to approximate the kernel function $\omega_{1-\alpha}(t)$ efficiently on the interval $[\Delta{t},\,T]$,
see \cite[Theorem 2.5]{Jiang2017Fast}.
\begin{lemma}\label{SOE}
For the given $\alpha\in(0,\,1)$, an absolute tolerance error $\epsilon\ll{1}$, a cut-off time $\Delta{t}>0$ and a finial time $T$, there exists a positive integer $N_{q}$, positive quadrature nodes $\theta^{\ell}$ and corresponding positive weights $\varpi^{\ell}\,(1\leq{\ell}\leq{N_{q}})$ such that
\begin{align}
\bigg|
\omega_{1-\alpha}(t)
-\sum_{\ell=1}^{N_{q}}\varpi^{\ell}e^{-\theta^{\ell}t}\bigg|\leq\epsilon,
\quad
\forall\,{t}\in[\Delta{t},\,T].\nonumber
\end{align}
\end{lemma}

The Caputo derivative \eqref{CaputoDef} is split into the sum of a history part (an integral over $[0,\,t_{n-1}]$) and a local part (an integral over $[t_{n-1},\,t_{n}]$) at the time $t_{n}$, see also \cite{Liao2018Unconditional}. Then, the local part will be approximated by linear interpolation directly, the history part can be evaluated via the SOE technique, that is,
\begin{align}\label{Fast-Approximation}
(\partial_{\tau}^{\alpha}v)(t_{n-\frac{1}{2}})
&\approx
\frac{1}{\tau_{n}}\int_{t_{n-1}}^{t_{n}}
\int_{0}^{t_{n-1}}v^{\prime}(s)
\sum_{\ell=1}^{N_{q}}
\varpi^{\ell}e^{-\theta^{\ell}(t-s)}\zd{s}\zd{t}\nonumber\\
&+\frac{1}{\tau_{n}}\int_{t_{n-1}}^{t_{n}}
\int_{t_{n-1}}^{t}\omega_{1-\alpha}(t-s)(\Pi_{1}v)^{\prime}(s)\zd{s}\zd{t}\nonumber\\
&=\bar{a}_{0}^{(n)}\triangledown_{\tau}v^{n}
+\frac{1}{\tau_{n}}\sum_{\ell=1}^{N_{q}}\varpi^{\ell}
\int_{t_{n-1}}^{t_{n}}
\int_{0}^{t_{n-1}}v^{\prime}(s)
e^{-\theta^{\ell}(t-t_{n-1})}
e^{-\theta^{\ell}(t_{n-1}-s)}\zd{s}\zd{t}\nonumber\\
&=\bar{a}_{0}^{(n)}\triangledown_{\tau}v^{n}
+\frac{1}{\tau_{n}}\sum_{\ell=1}^{N_{q}}\varpi^{\ell}
b^{(n,l)}
\mathcal{H}^{\ell}(t_{n-1}),\quad n\geq{1},
\end{align}
where
\[
\mathcal{H}^{\ell}(t_{k})
:=\int_{0}^{t_{k}}e^{-\theta^{\ell}(t_{k}-s)}v^{\prime}(s)\zd{s}\quad
\text{with}\quad\mathcal{H}^{\ell}(t_{0})=0,\quad
b^{(k,l)}:=\int_{t_{k-1}}^{t_{k}}e^{-\theta^{\ell}(t-t_{k-1})}\zd{t}.
\]
By utilizing the linear interpolation and a recursive formula, we can approximate $\mathcal{H}^{\ell}(t_{k})$ by
\begin{align}\label{HistoryPart}
\mathcal{H}^{\ell}(t_{k})
&\approx\int_{0}^{t_{k-1}}e^{-\theta^{\ell}(t_{k}-s)}v^{\prime}(s)\zd{s}
+\int_{t_{k-1}}^{t_{k}}e^{-\theta^{\ell}(t_{k}-s)}
\frac{\triangledown_{\tau}v^{k}}{\tau_{k}}\zd{s}\nonumber\\
&=e^{-\theta^{\ell}\tau_{k}}\mathcal{H}^{\ell}(t_{k-1})
+c^{(k,l)}\triangledown_{\tau}v^{k},
\end{align}
where the positive coefficients
\[
c^{(k,l)}:=\frac{1}{\tau_{k}}\int_{t_{k-1}}^{t_{k}}e^{-\theta^{\ell}(t_{k}-s)}\zd{s},\quad k\geq{1}.
\]
Having taken this excursion through {\eqref{Fast-Approximation}}-{\eqref{HistoryPart}},
 we arrive at the fast L1$^{+}$ formula
\begin{align}\label{Fast-NewL1-Formula}
(\partial_{f}^{\alpha}v)^{n-\frac{1}{2}}
=\bar{a}_{0}^{(n)}\triangledown_{\tau}v^{n}
+\frac{1}{\tau_{n}}\sum_{\ell=1}^{N_{q}}\varpi^{\ell}
b^{(n,l)}\mathcal{H}^{\ell}(t_{n-1}),
\end{align}
in which $\mathcal{H}^{\ell}(t_{k})$ is computed by using the recursive relationship
\begin{align}\label{HistoryRecursive}
\mathcal{H}^{\ell}(t_{k})
=e^{-\theta^{\ell}\tau_{k}}\mathcal{H}^{\ell}(t_{k-1})
+c^{(k,l)}\triangledown_{\tau}v^{k},\quad k\geq{1},\quad 1\leq\ell\leq{N_{q}}.
\end{align}
\subsection{Adaptive time-stepping strategy}
In the previous sections, we have proved that the numerical schemes are unconditionally energy stable which implies large time steps are allowed.
Indeed, in simulating the phase field problems such as the coarsening dynamics problems discussed in Example \ref{Coarsening-Dynamics},
adaptive time-stepping strategy is necessary to efficiently resolve widely varying time scales
and to significantly reduce the computational cost.
Roughly speaking, the adaptive time steps can be selected by using an accuracy criterion example as \cite{Gomez2011Provably}, or
the time evolution of the total energy such as \cite{Qiao2011An}.
We focus on the former and update the time step size by using the formula
\begin{align*}
\tau_{ada}\bra{e,\tau}
=\rho\bra{\frac{tol}{e}}^{\frac{1}{2}}\tau,
\end{align*}
where $\rho$ is a default safety coefficient, $tol$ is a reference tolerance, and $e$ is the relative error at each time level.
The details of the adaptive time steps strategy are presented in Algorithm \ref{Adaptive-Time-Step-Strategy}.
Here, the first-order SAV and  second-order SAV schemes refer to the backward Euler-L1 method
and Crank-Nicolson SAV method proposed in this article, respectively.

\begin{algorithm}
\caption{Adaptive time-stepping strategy}
\label{Adaptive-Time-Step-Strategy}
\begin{algorithmic}[1]
\Require{Given $\phi^{n}$ and time step $\tau_{n}$}
\State Compute $\phi_{1}^{n+1}$ by using first-order SAV scheme with time step $\tau_{n}$.
\State Compute $\phi_{2}^{n+1}$ by using second-order SAV scheme with time step $\tau_{n}$.
\State Calculate $e_{n+1}=\|\phi_{2}^{n+1}-\phi_{1}^{n+1}\|/\|\phi_{2}^{n+1}\|$.
\If {$e_{n}<tol$ or $\tau_{n}\leq{\tau_{\min}}}$
\If {$e_{n}<tol$}
\State Update time-step size $\tau_{n+1}\leftarrow\min\{\max\{\tau_{\min},\tau_{ada}\},\tau_{\max}\}$.
\Else
\State Update time-step size $\tau_{n+1}\leftarrow\tau_{\min}$.
\EndIf
\Else
\State Recalculate with time-step size  $\tau_{n}\leftarrow\min\{\max\{\tau_{\min},\tau_{ada}\},\tau_{\max}\}$.
\State Goto 1
\EndIf
\end{algorithmic}
\end{algorithm}

\subsection{Numerical examples}
The CN-SAV methods \eqref{Slope-CN-SAV-1}-\eqref{Slope-CN-SAV-2} and \eqref{No-CN-SAV-1}-\eqref{No-CN-SAV-2} are examined in this section for the time-fractional MBE model \eqref{Problem-2}.
Specially, the model fractional ODE problem
$\partial_{t}^{\alpha}u(t)=f(t)$ is used in Example \ref{Accuracy-Test-New-L1} to test the accuracy
of L1$^{+}$ formula \eqref{New-L1-Formula}.
The time interval $[0,T]$ is always divided into two parts $[0, T_{0}]$ and $[T_{0}, T]$ with total $N$ subintervals.
We will
take $T_0=\min\{1/\gamma,T\}$, and
apply the graded grid $t_{k}=T_{0}(k/N_0)^{\gamma}$
in $[0,T_{0}]$ to resolve the initial singularity.
In the remainder interval $[T_{0},T]$,
we put $N_1:=N-N_0$ cells with random time-steps
$$\tau_{N_{0}+k}=\frac{(T-T_{0})\epsilon_{k}}{\sum_{k=1}^{N_1}\epsilon_{k}}\quad\text{for $1\leq k\leq N_1$}$$
where $\epsilon_{k}\in(0,1)$ are the random numbers.

The time accuracy of the proposed methods is mainly focused on, so for simplicity,
and the Fourier pseudo-spectral method is always applied to approximate
the space variables using the same spacing in each spatial direction.
To examine the CN-SAV schemes,
the maximum norm error $e(N):=\max_{1\leq{n}\leq{N}}\|U^{n}-u^{n}\|_{\infty}$ is recorded in each run,
and the experimental convergence order in time is computed by
$$\text{Order}:=\frac{\log\bra{e(N)/e(2N)}}{\log\bra{\tau(N)/\tau(2N)}}$$
where $\tau(N)$ denotes the maximum time-step size for total $N$ subintervals

\begin{example}\label{Accuracy-Test-New-L1}
Consider the model fractional ODE problem
$\partial_{t}^{\alpha}u(t)=f(t)$ with an exact solution $u=\omega_{1+\sigma}(t)$,
where the parameter $\sigma$ determines the initial regularity of $u$.
The accuracy of L1$^{+}$ formula \eqref{New-L1-Formula} is examined carefully
using the following three scenarios:
\begin{enumerate}
\item[(a)] $\sigma>{2}$ and $\alpha=0.1,\,0.5$ and $0.9$ on uniform mesh, see Table \ref{New-L1-Error-1} in subsection 3.1;
\item[(b)] $\sigma<{2}$ and $\alpha=0.1,\,0.5$ and $0.9$ using the graded parameter $\gamma=1$, see Table \ref{New-L1-Error-2};
\item[(c)] $\sigma<{2}$ and $\alpha=0.3$ and $0.7$ using different graded parameters $\gamma\geq1$, see Tables \ref{New-L1-Error-3}-\ref{New-L1-Error-4}.
\end{enumerate}
\end{example}

\begin{table}[htb!]
\begin{center}
\caption{Numerical accuracy of L1$^{+}$ formula \eqref{New-L1-Formula} with $\sigma=0.8,\gamma=1$}\label{New-L1-Error-2}
\def\temptablewidth{1.0\textwidth}
{\rule{\temptablewidth}{0.5pt}}
\begin{tabular*}{\temptablewidth}{@{\extracolsep{\fill}}cccccccc}
\multirow{2}{*}{$N$} &\multirow{2}{*}{$\tau$} &\multicolumn{2}{c}{$\alpha=0.1$} &\multicolumn{2}{c}{$\alpha=0.5$} &\multicolumn{2}{c}{$\alpha=0.9$}\\
               \cline{3-4}       \cline{5-6}       \cline{7-8}
        &       &$e(N)$  &Order     &$e(N)$ &Order     &$e(N)$  &Order \\
\midrule
  64     &2.24e-02     &5.90e-03	 &$-$      &4.65e-03	&$-$     &1.10e-03	 &$-$\\
  128    &1.16e-02     &3.39e-03	 &0.84      &2.67e-03	&0.84    &6.29e-04	 &0.84\\
  256    &5.87e-03     &1.95e-03	 &0.82      &1.53e-03	&0.82    &3.61e-04	 &0.82\\
  512    &2.88e-03     &1.12e-03	 &0.78      &8.80e-04	&0.78    &2.07e-04	 &0.78\\
\end{tabular*}
{\rule{\temptablewidth}{0.5pt}}
\end{center}
\end{table}

\begin{table}[htb!]
\begin{center}
\caption{Numerical accuracy of L1$^{+}$ formula \eqref{New-L1-Formula} with $\alpha=0.3,\,\sigma=0.5$}\label{New-L1-Error-3} \vspace*{0.3pt}
\def\temptablewidth{1.0\textwidth}
{\rule{\temptablewidth}{0.5pt}}
\begin{tabular*}{\temptablewidth}{@{\extracolsep{\fill}}cccccccccc}
\multirow{2}{*}{$N$} &\multirow{2}{*}{$\tau$} &\multicolumn{2}{c}{$\gamma=2$} &\multirow{2}{*}{$\tau$} &\multicolumn{2}{c}{$\gamma=4$} &\multirow{2}{*}{$\tau$}&\multicolumn{2}{c}{$\gamma=5$} \\
             \cline{3-4}          \cline{6-7}         \cline{9-10}
         &          &$e(N)$   &Order &         &$e(N)$   &Order &         &$e(N)$    &Order\\
\midrule
  64     &2.78e-02  &6.47e-03 &$-$   &2.68e-02 &1.21e-03 &$-$   &2.62e-02 &1.68e-03  &$-$\\
  128    &1.53e-02  &3.23e-03 &1.15  &1.38e-02 &3.70e-04 &1.77  &1.30e-02 &5.04e-04	 &1.71\\
  256    &7.37e-03  &1.63e-03 &0.94  &6.57e-03 &8.15e-05 &2.04  &6.76e-03 &1.55e-04	 &1.81\\
  512    &3.64e-03  &8.20e-04 &0.97  &3.51e-03 &2.46e-05 &1.91  &3.27e-03 &3.77e-05	 &1.95\\
\end{tabular*}
{\rule{\temptablewidth}{0.5pt}}
\end{center}
\end{table}	

\begin{table}[htb!]
\begin{center}
\caption{Numerical accuracy of L1$^{+}$ formula \eqref{New-L1-Formula} with $\alpha=0.7,\,\sigma=0.5$}\label{New-L1-Error-4} \vspace*{0.3pt}
\def\temptablewidth{1.0\textwidth}
{\rule{\temptablewidth}{0.5pt}}
\begin{tabular*}{\temptablewidth}{@{\extracolsep{\fill}}cccccccccc}
\multirow{2}{*}{$N$} &\multirow{2}{*}{$\tau$} &\multicolumn{2}{c}{$\gamma=2$} &\multirow{2}{*}{$\tau$} &\multicolumn{2}{c}{$\gamma=4$} &\multirow{2}{*}{$\tau$}&\multicolumn{2}{c}{$\gamma=5$} \\
             \cline{3-4}          \cline{6-7}         \cline{9-10}
         &          &$e(N)$   &Order &         &$e(N)$   &Order &         &$e(N)$    &Order\\
\midrule
  64     &2.76e-02  &8.40e-04 &$-$   &2.68e-02 &3.46e-04 &$-$   &2.71e-02 &6.66e-04  &$-$\\
  128    &1.49e-02  &4.20e-04 &1.12  &1.37e-02 &9.47e-05 &1.93  &1.30e-02 &1.80e-04	 &1.78\\
  256    &7.34e-03  &2.12e-04 &0.97  &6.90e-03 &2.47e-05 &1.95  &6.76e-03 &4.91e-05	 &1.98\\
  512    &3.65e-03  &1.07e-04 &0.99  &3.47e-03 &6.44e-06 &1.96  &3.26e-03 &1.19e-05	 &1.94\\
\end{tabular*}
{\rule{\temptablewidth}{0.5pt}}
\end{center}
\end{table}

Table \ref{New-L1-Error-1} lists the case (a) having smooth solution, while Tables \ref{New-L1-Error-2}-\ref{New-L1-Error-4}
record the two cases (b)-(c) having non-smooth solutions.
From these numerical results in Tables \ref{New-L1-Error-1}-\ref{New-L1-Error-4}, one sees that it is accurate of $O(\tau^{\max\{\gamma\sigma,2\}})$
via the following observations:
(i) The numerical accuracy is independent of the fractional order $\alpha\in(0,1)$,
and it is second-order accurate for smooth solutions with $\sigma\geq{2}$;
(ii) On the uniform mesh, the numerical accuracy degenerates to $O(\tau^{\sigma})$ when the regularity parameter $\sigma\in(0,2)$;
(iii) When the solution is non-smooth, the numerical accuracy reaches $O(\tau^{\gamma\sigma})$ by the graded mesh,
and the second-order accuracy would be recovered by choosing $\gamma\geq{2/\sigma}$.

\begin{example}\label{Accuracy-Test-CN-SAV}
To examine the temporal accuracy of our CN-SAV schemes, consider the time-fractional MBE model
$\partial_{t}^{\alpha}\phi
=-M\bra{\varepsilon^{2}\Delta^{2}\phi+f\bra{\nabla\phi}}
+g(\mathbf{x},t)$
for $\mathbf{x}\in(0,2\pi)^{2}$ and $0<t<1$
such that it has an exact solution $\phi=\omega_{1+\sigma}(t)\sin(x)\sin(y)$.
\end{example}

\begin{table}[htb!]
\begin{center}
\caption{Numerical accuracy of CN-SAV scheme \eqref{Slope-CN-SAV-1}-\eqref{Slope-CN-SAV-2}  with $\alpha=0.8,\,\sigma=0.4$}\label{Slope-CN-SAV-Error} \vspace*{0.3pt}
\def\temptablewidth{1.0\textwidth}
{\rule{\temptablewidth}{0.5pt}}
\begin{tabular*}{\temptablewidth}{@{\extracolsep{\fill}}cccccccccc}
\multirow{2}{*}{$N$} &\multirow{2}{*}{$\tau$} &\multicolumn{2}{c}{$\gamma=3$} &\multirow{2}{*}{$\tau$} &\multicolumn{2}{c}{$\gamma=5$} &\multirow{2}{*}{$\tau$}&\multicolumn{2}{c}{$\gamma=6$} \\
             \cline{3-4}          \cline{6-7}         \cline{9-10}
             &     &$e(N)$   &Order    &      &$e(N)$   &Order    &      &$e(N)$    &Order\\
\midrule
  64     &3.68e-02 &1.78e-03 &$-$  &3.92e-02  &5.05e-04 &$-$  &4.04e-02	 &4.63e-04  &$-$\\
  128    &1.76e-02 &7.87e-04 &1.11 &2.10e-02  &1.33e-04 &2.13 &2.25e-02	 &1.19e-04	&2.33\\
  256    &9.75e-03 &3.43e-04 &1.40 &1.07e-02  &3.46e-05 &2.00 &1.07e-02  &2.96e-05  &1.86\\
  512    &4.82e-03 &1.49e-04 &1.18 &5.17e-03  &8.80e-06 &1.89 &5.33e-03	 &7.39e-06	&2.00\\
\end{tabular*}
{\rule{\temptablewidth}{0.5pt}}
\end{center}
\end{table}

\begin{table}[htb!]
\begin{center}
\caption{Numerical accuracy of CN-SAV scheme \eqref{No-CN-SAV-1}-\eqref{No-CN-SAV-2} with $\alpha=0.8,\,\sigma=0.4$}\label{No-CN-SAV-Error} \vspace*{0.3pt}
\def\temptablewidth{1.0\textwidth}
{\rule{\temptablewidth}{0.5pt}}
\begin{tabular*}{\temptablewidth}{@{\extracolsep{\fill}}cccccccccc}
\multirow{2}{*}{$N$} &\multirow{2}{*}{$\tau$} &\multicolumn{2}{c}{$\gamma=3$} &\multirow{2}{*}{$\tau$} &\multicolumn{2}{c}{$\gamma=5$} &\multirow{2}{*}{$\tau$}&\multicolumn{2}{c}{$\gamma=6$} \\
             \cline{3-4}          \cline{6-7}         \cline{9-10}
  &      &$e(N)$   &Order &         &$e(N)$   &Order &         &$e(N)$    &Order\\
\midrule
  64     &3.96e-02 &1.78e-03 &$-$   &4.38e-02 &5.04e-04 &$-$  &4.26e-02 &4.63e-04 &$-$\\
  128    &2.02e-02 &7.87e-04 &1.21  &2.16e-02 &1.33e-04 &1.88 &2.23e-02 &1.17e-04 &2.14\\
  256    &9.74e-03 &3.43e-04 &1.14  &1.06e-02 &3.46e-05 &1.89 &1.07e-02	&2.87e-05 &1.90	\\
  512    &4.96e-03 &1.49e-04 &1.23  &5.37e-03 &8.80e-06 &2.02 &5.51e-03	&7.52e-06 &2.03\\
\end{tabular*}
{\rule{\temptablewidth}{0.5pt}}
\end{center}
\end{table}

The parameters are taken as $M=0.1,\beta=1,C_{0}=1,\varepsilon^{2}=0.5$, and the space is discretized by $128\times{128}$ meshes.
We run the CN-SAV schemes \eqref{Slope-CN-SAV-1}-\eqref{Slope-CN-SAV-2} and \eqref{No-CN-SAV-1}-\eqref{No-CN-SAV-2} by setting a variety of regularity parameters.
Numerical results are tabulated in Tables \ref{Slope-CN-SAV-Error}-\ref{No-CN-SAV-Error}, respectively.
Tables \ref{Slope-CN-SAV-Error} and \ref{No-CN-SAV-Error} show
the numerical results in the worse case of $\sigma<\alpha$.
It seen that it is accurate of order $O(\tau^{\gamma\sigma})$ on the graded mesh,
and the second-order accuracy can be recovered by taking $\gamma\geq\gamma_{\mathrm{opt}}=2/\sigma$.
The computational results suggest that it is convergent of $O(\tau^{\min\{\gamma\sigma,2\}})$ in time
although no theoretical proof is available up to now.

\begin{example}\label{Time-Adaptive-Test-CN-SAV}
We examine here the performance of the adaptive time-stepping method for the CN-SAV schemes \eqref{Slope-CN-SAV-1}-\eqref{Slope-CN-SAV-2} and \eqref{No-CN-SAV-1}-\eqref{No-CN-SAV-2} with parameter value $\epsilon^{2}=0.1$.
To this end, we carry out a standard benchmark problem using initial condition
\begin{align}\label{Given-Initial-Condition}
\phi\bra{\mathbf{x},0}
=0.1\bra{\sin(3x)\sin(2y)+\sin(5x)\sin(5y)}.
\end{align}
\end{example}

In what follows, if not explicitly specified, the default values of parameters are given as
$M=1,\beta=4,\mathbf{x}\in(0,2\pi)^{2},\rho=0.9, tol=10^{-3}$, and the space is discretized by $128\times{128}$ meshes.
To quantify the deviation of the height function, define the roughness function $W(t)$:
\begin{align}
W(t)
=\sqrt{\frac{1}{\abs{\Omega}}
\int_{\Omega}\bra{\phi\bra{\mathbf{x},t}
-\bar{\phi}(t)}^{2}\zd{\mathbf{x}}},
\end{align}
where $\bar{\phi}(t)=\frac{1}{\abs{\Omega}}
\int_{\Omega}\phi\bra{\mathbf{x},t}\zd{\mathbf{x}}$ is the average. The function $W(t)$ characterizes the mean size of the network cell.

\begin{figure}[htb!]
\centering
\includegraphics[width=3.0in,height=2.0in]{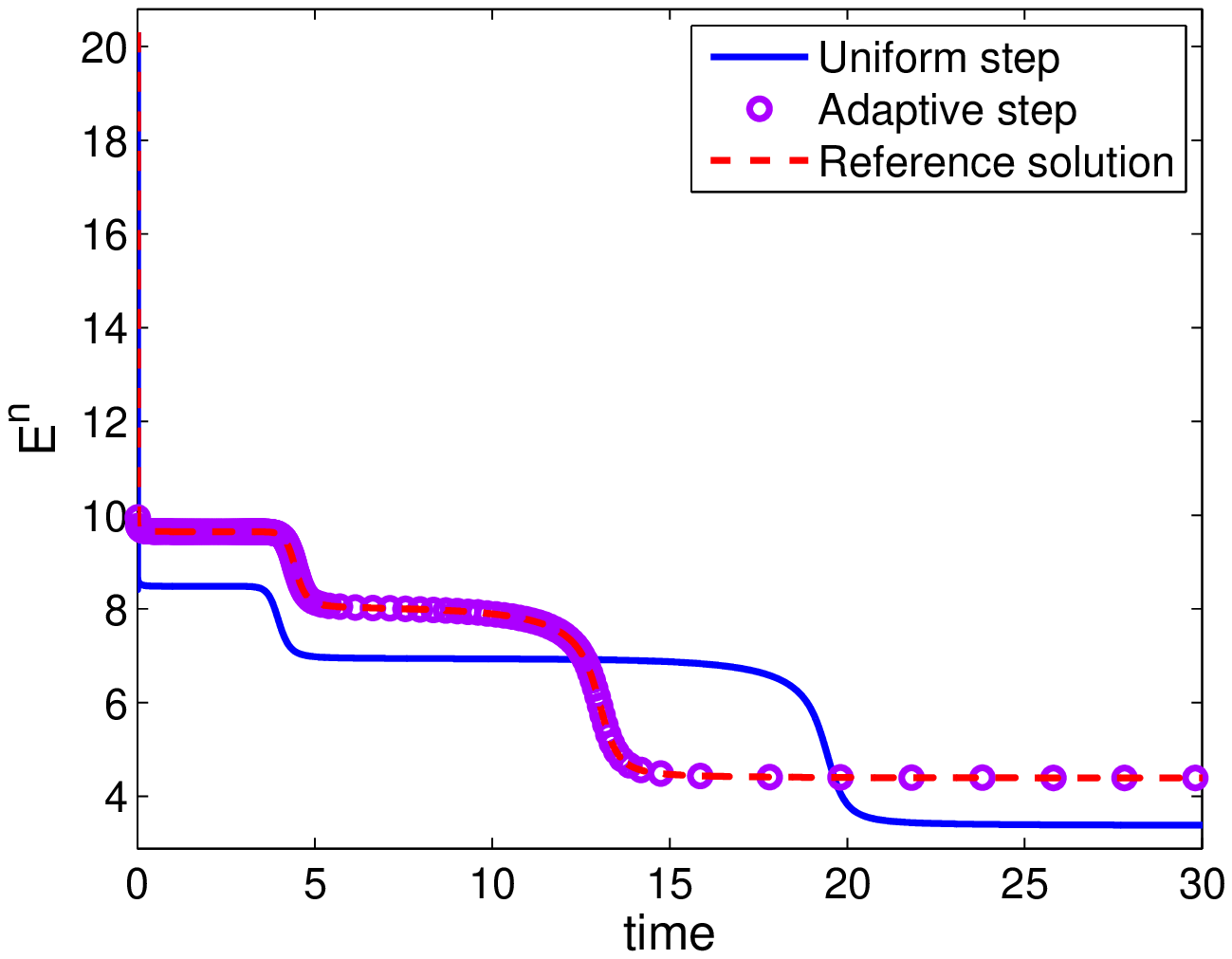}
\includegraphics[width=3.0in,height=2.0in]{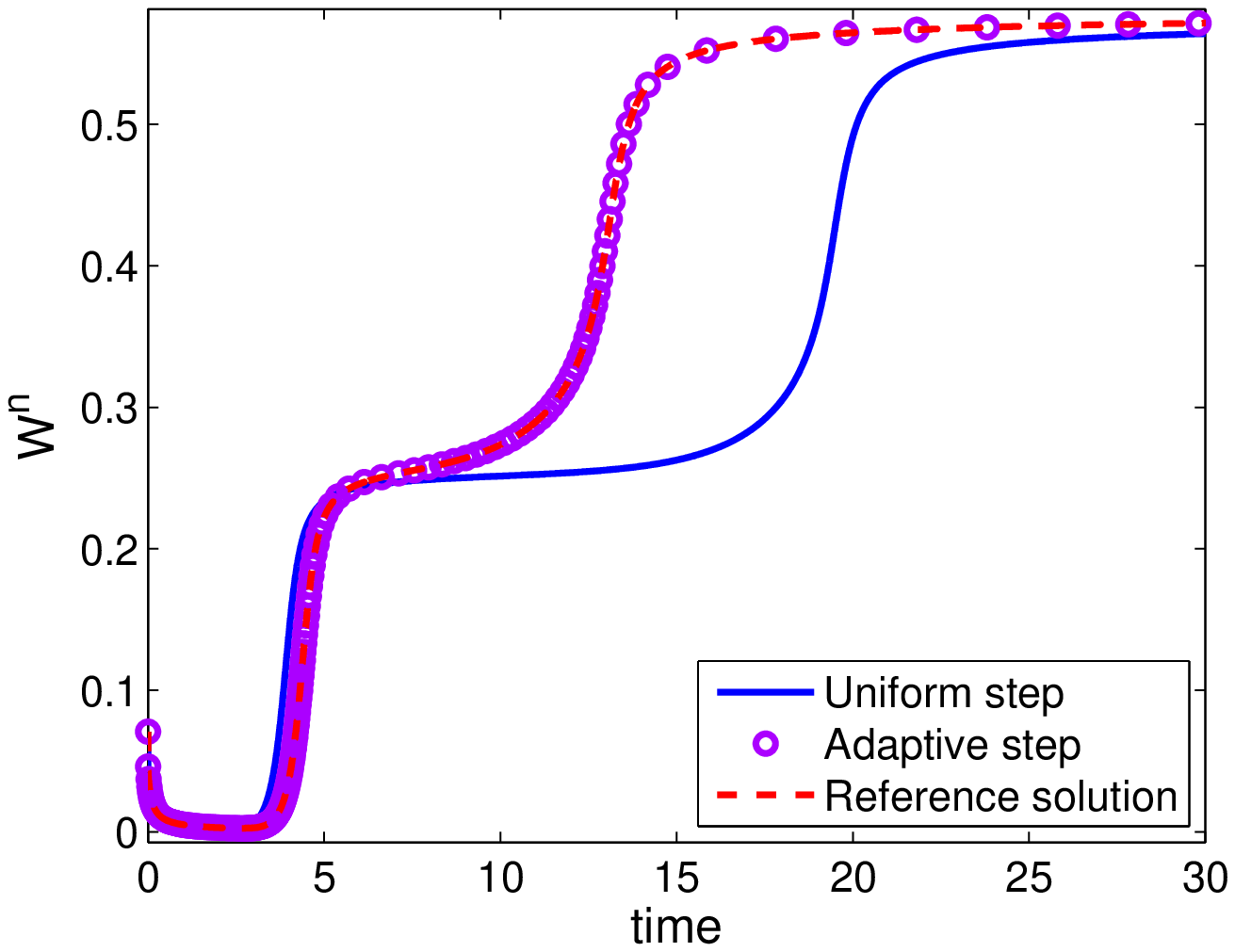}\\
\includegraphics[width=3.0in,height=2.0in]{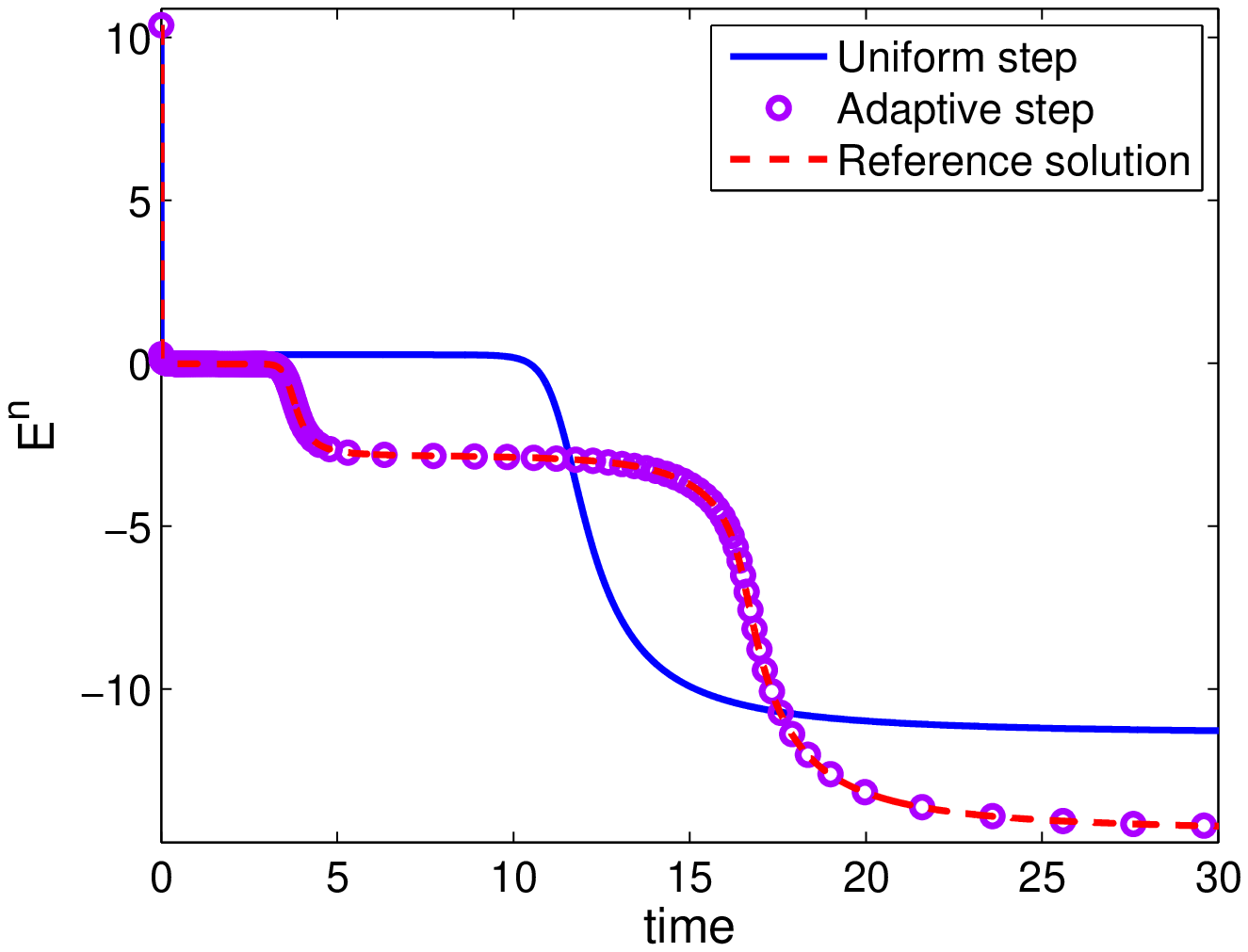}
\includegraphics[width=3.0in,height=2.0in]{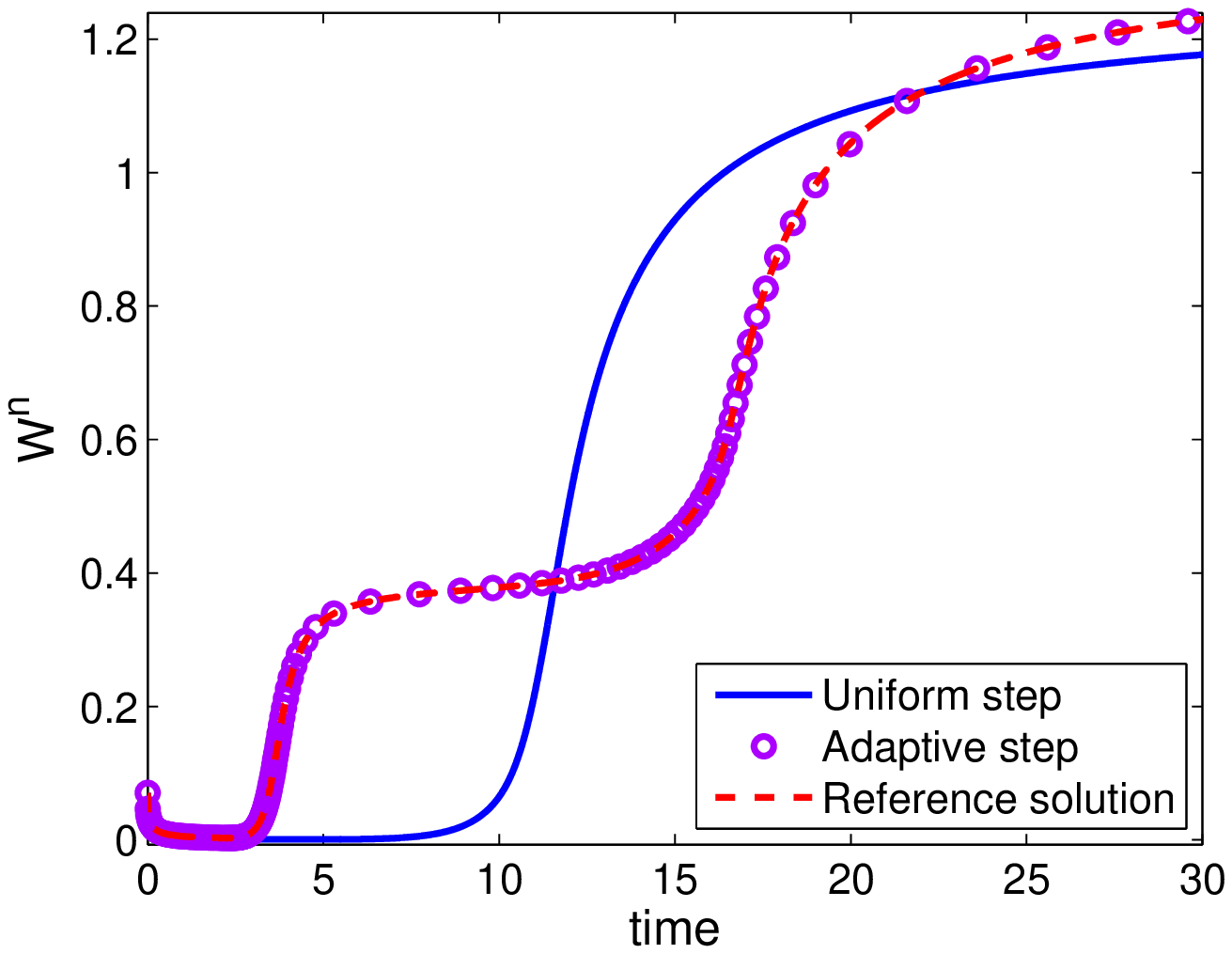}\\
\caption{Evolutions of energy and roughness (from left to right) of the Slope-Model and the No-Slope-Model (from top to bottom)
for $\alpha=0.7$ using different time steps strategies.}
\label{Slope-Model-Adaptive-Test}
\end{figure}

\begin{figure}[htb!]
\centering
\includegraphics[width=3.0in,height=2.0in]{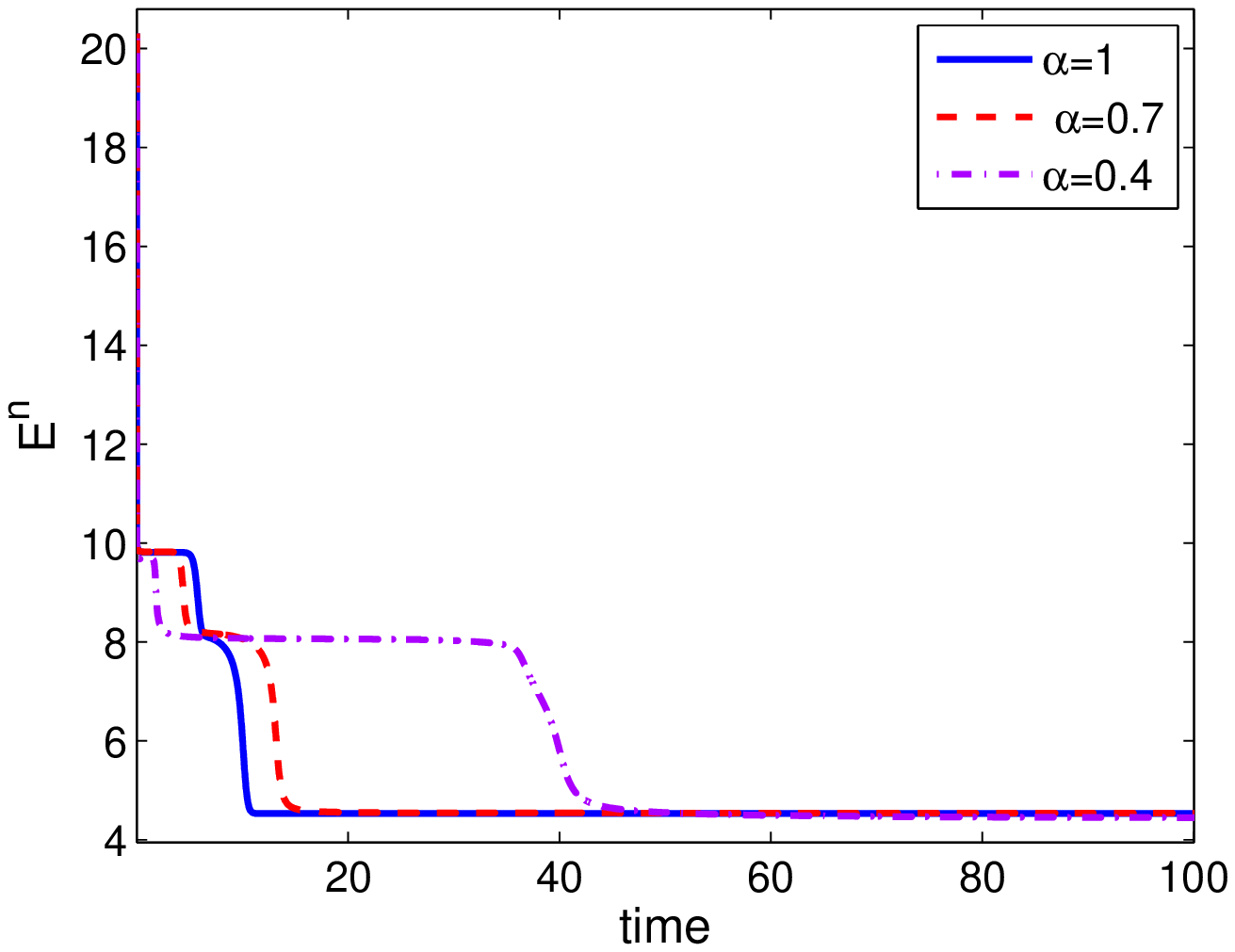}
\includegraphics[width=3.0in,height=2.0in]{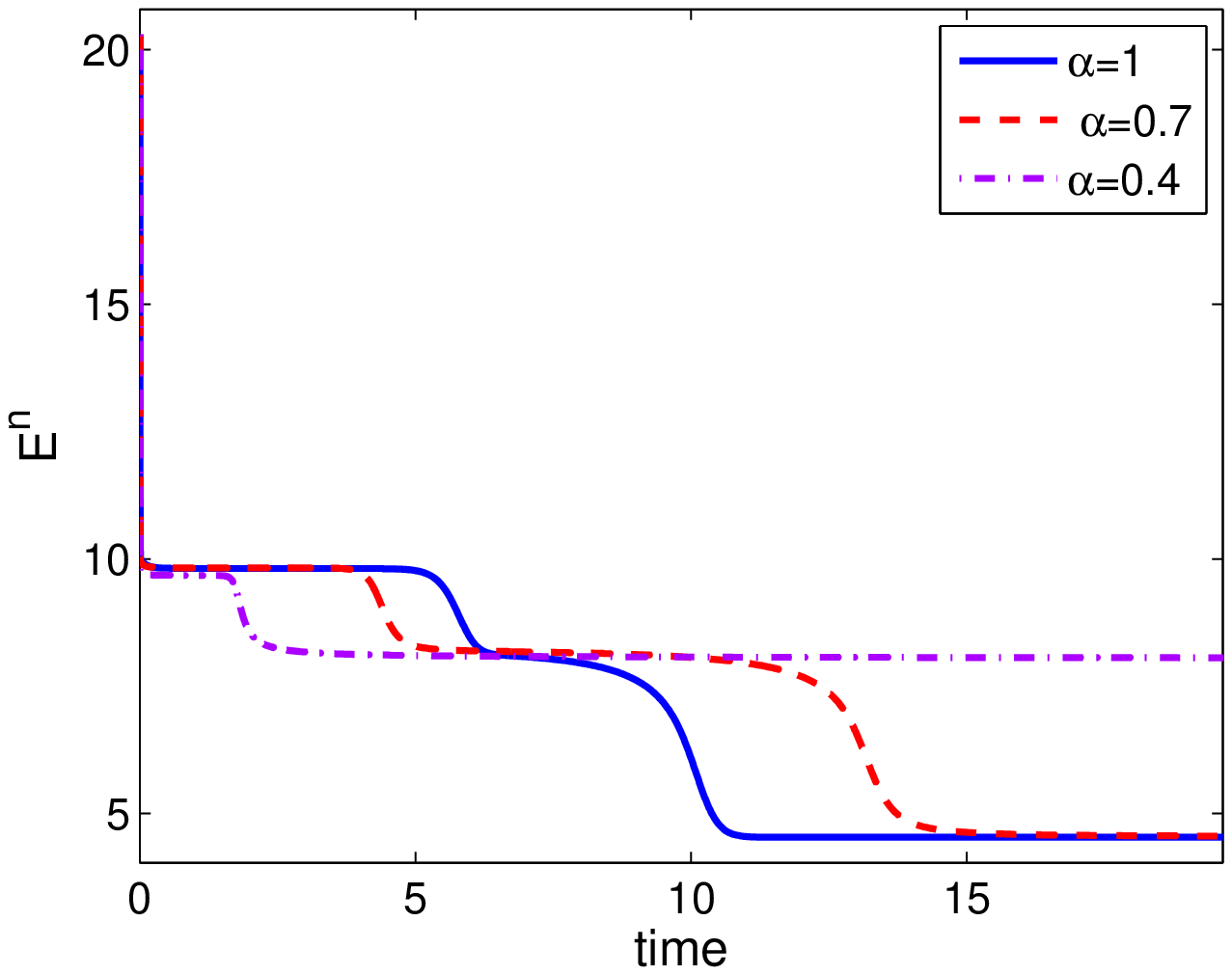}\\
\includegraphics[width=3.0in,height=2.0in]{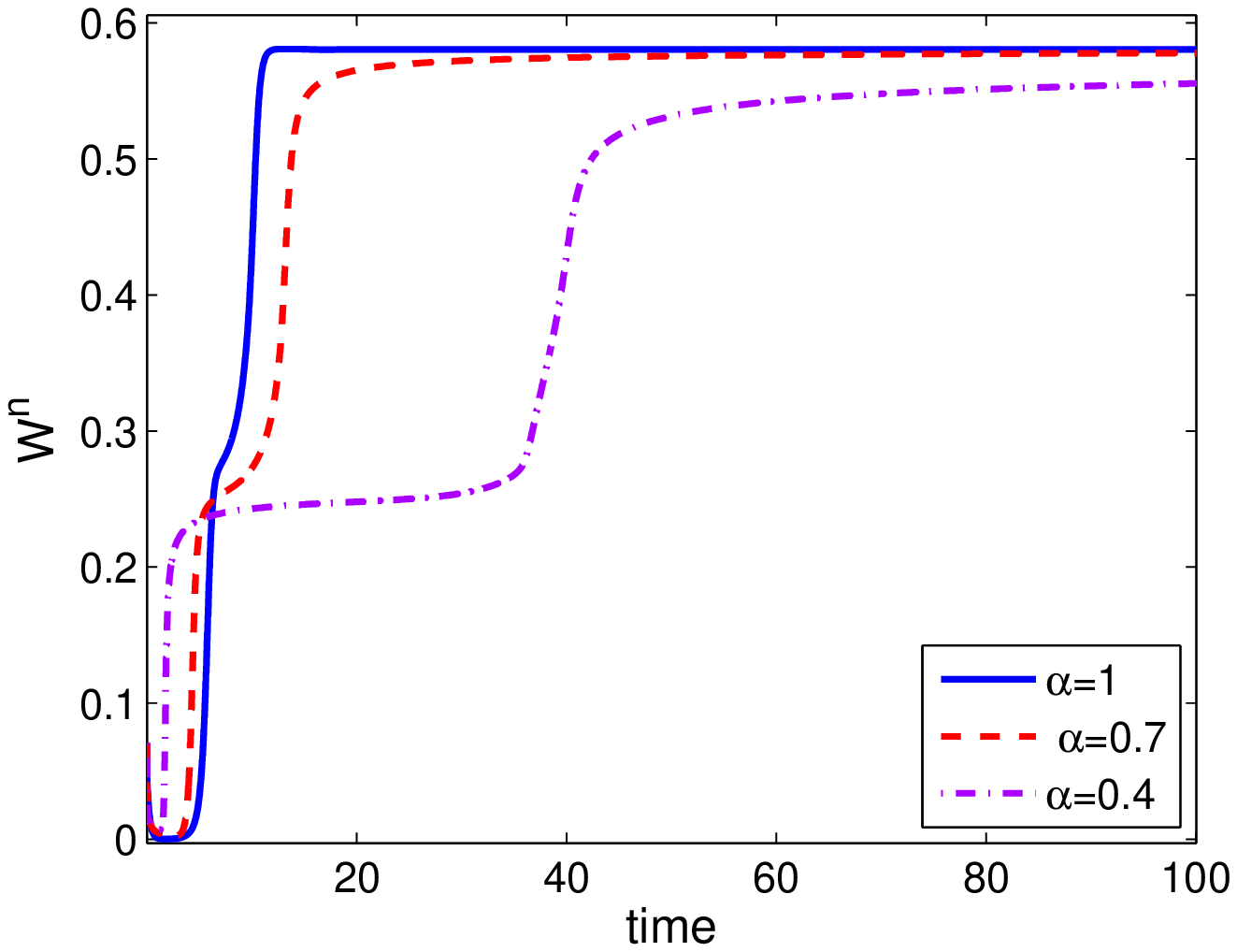}
\includegraphics[width=3.0in,height=2.0in]{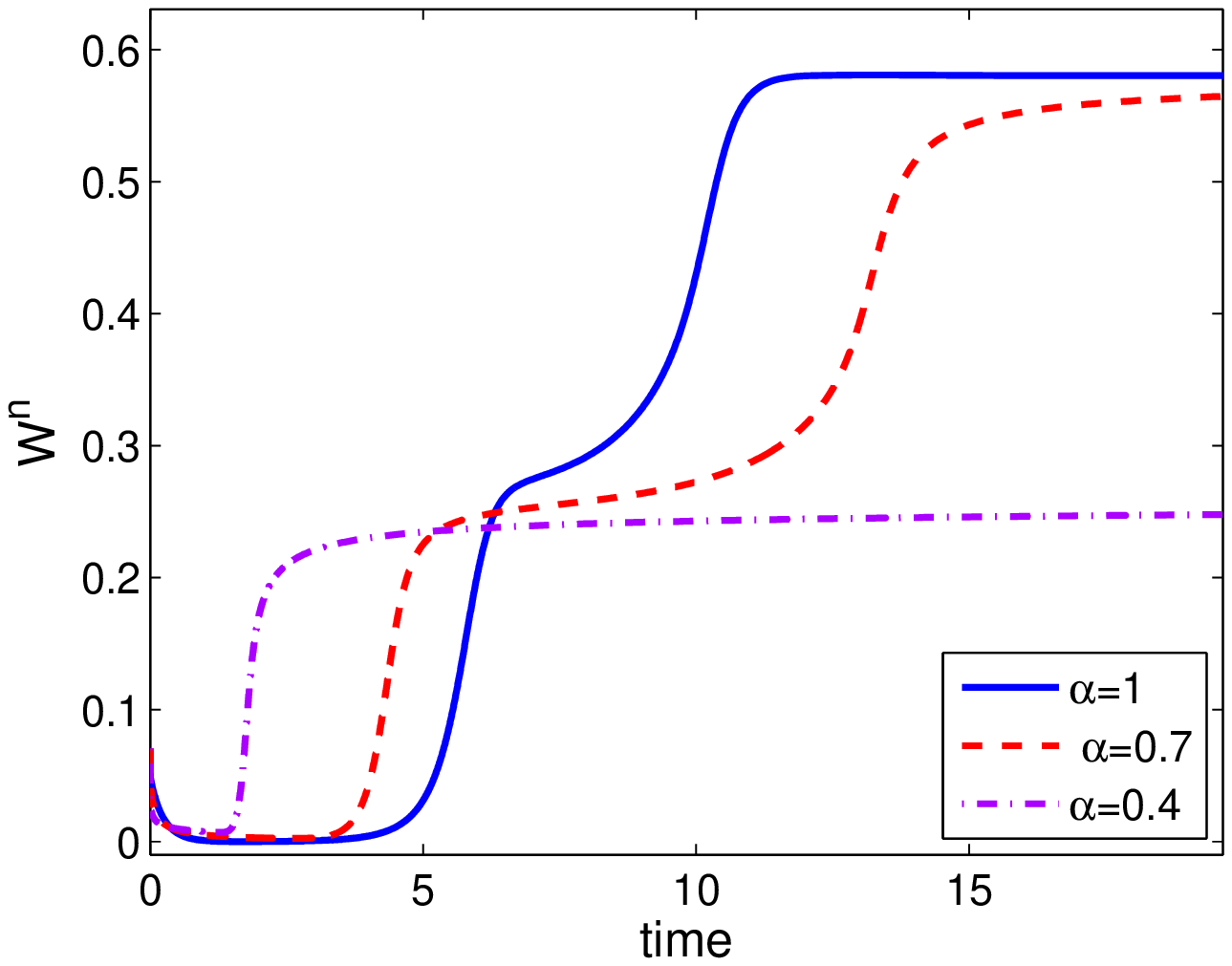}\\
\caption{Evolutions of energy and roughness of the Slope-Model (from top to bottom) for $t\in[0,100]$ and $t\in[0,20]$ (from left to right) with $\alpha=0.4,0.7$ and $1$, respectively.}
\label{Slope-Model-Energy-Roughness}
\end{figure}

For a fixed fractional index $\alpha=0.7$ and final time $T=30$, we first apply a constant time step $\tau=10^{-3}$, i.e., $N=30000$, to compute the solution.
Recall that the intrinsically initial singularity of solution that presented in Figure \ref{Initial-Singularity-alpha-04} early, the numerical results suggest the time mesh should be refined near the initial time.
As a consequence, we could obtain the reference solution where the parameter values $\gamma=3,T_{0}=0.01,N_{0}=30$ are applied in the cell $[0,T_{0}]$ and the uniform mesh is used over the remainder with the total numbers $N=30000$ as before.
For the adaptive time-stepping technique, taking the analogous numerical strategy in the initial time and choosing parameter $\tau_{\min}=\tau_{N_{0}}=10^{-3}$
and $\tau_{\max}=10^{-1}$ for the remainder.

\begin{figure}[htb!]
\centering
\includegraphics[width=3.0in,height=2.0in]{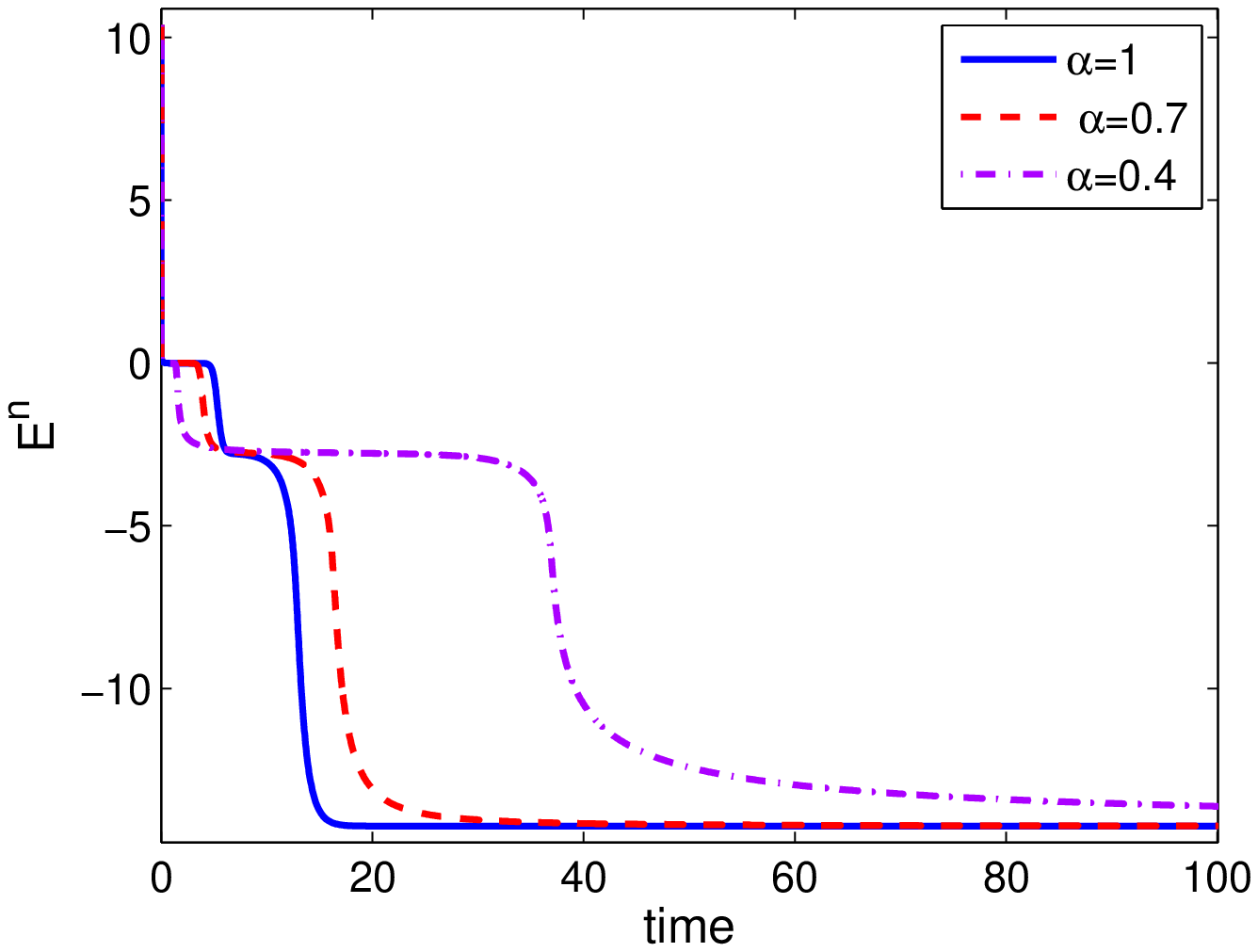}
\includegraphics[width=3.0in,height=2.0in]{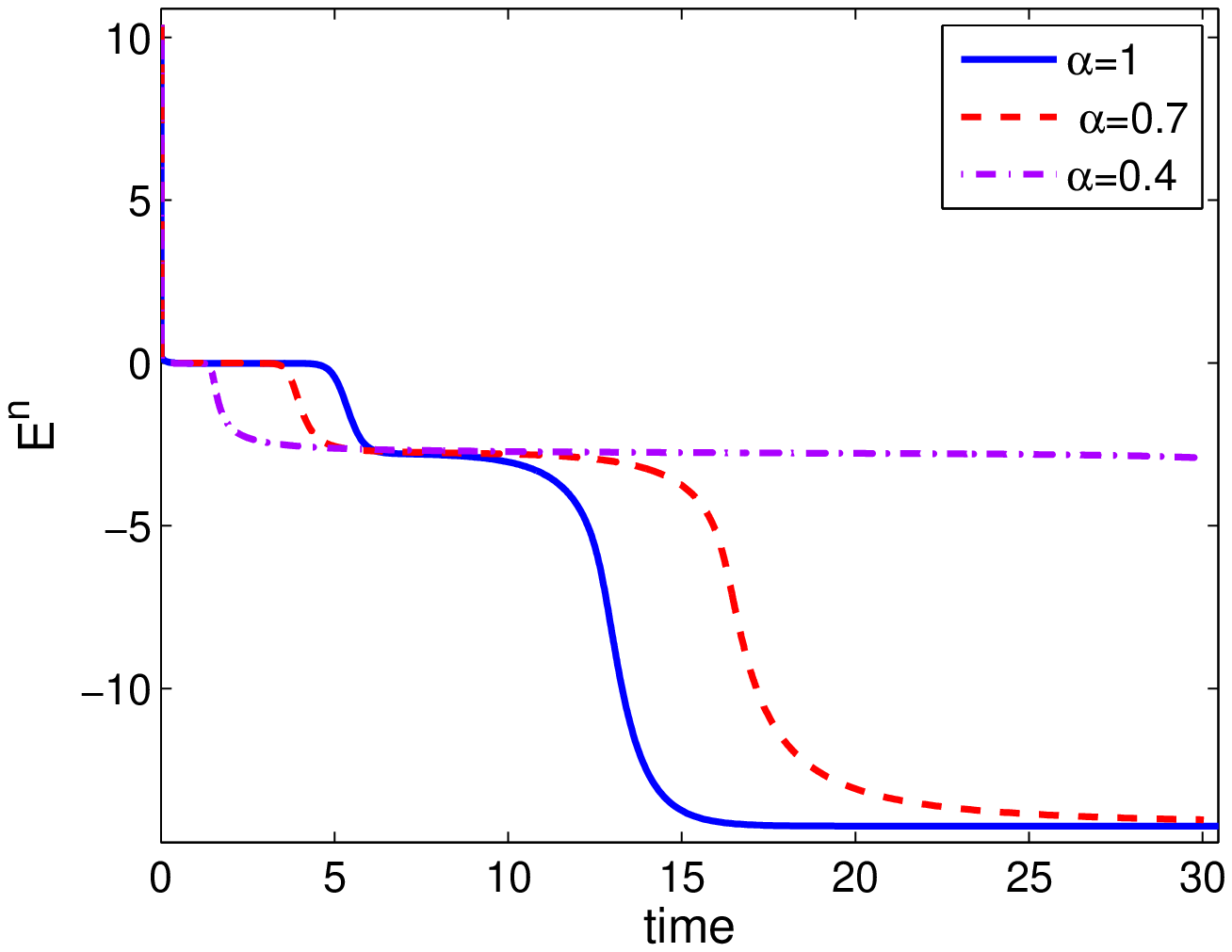}\\
\includegraphics[width=3.0in,height=2.0in]{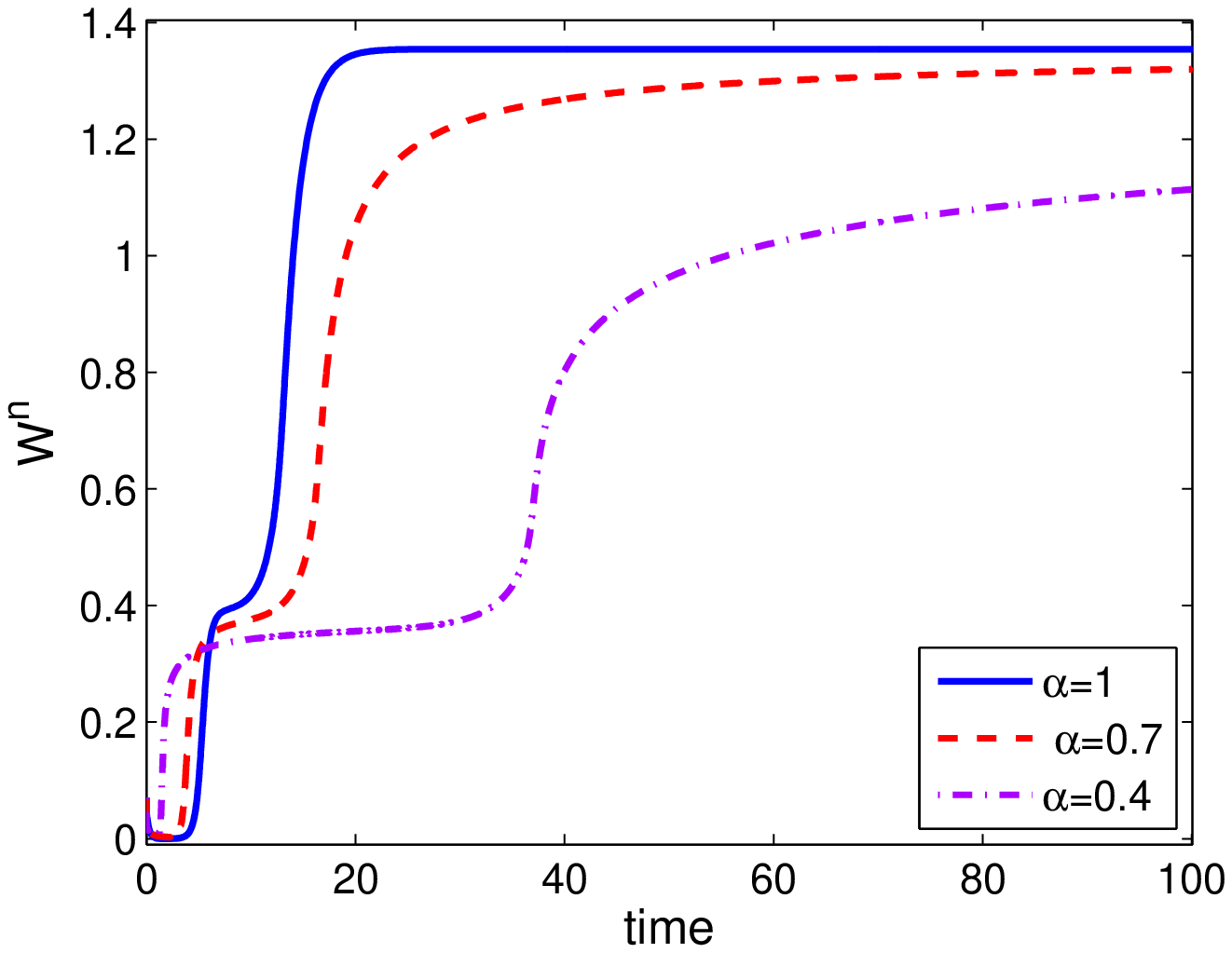}
\includegraphics[width=3.0in,height=2.0in]{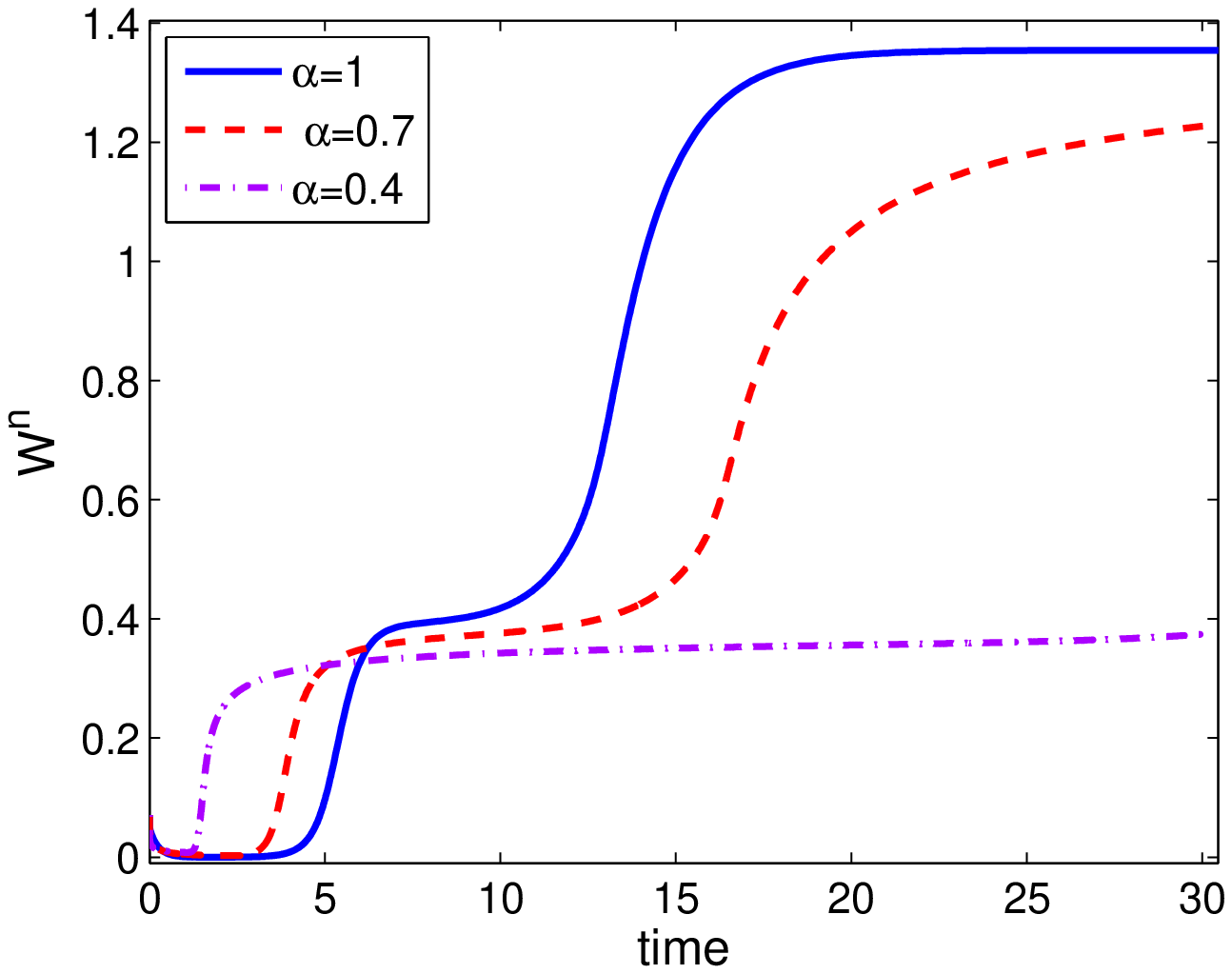}\\
\caption{Evolutions of energy and roughness of the No-Slope-Model (from top to bottom) for $t\in[0,100]$ and $t\in[0,30]$ (from left to right) with $\alpha=0.4,0.7$ and $1$, respectively.}
\label{No-Slope-Model-Energy-Roughness}
\end{figure}

\begin{figure}[htb!]
\centering
\includegraphics[width=3.0in,height=2.0in]{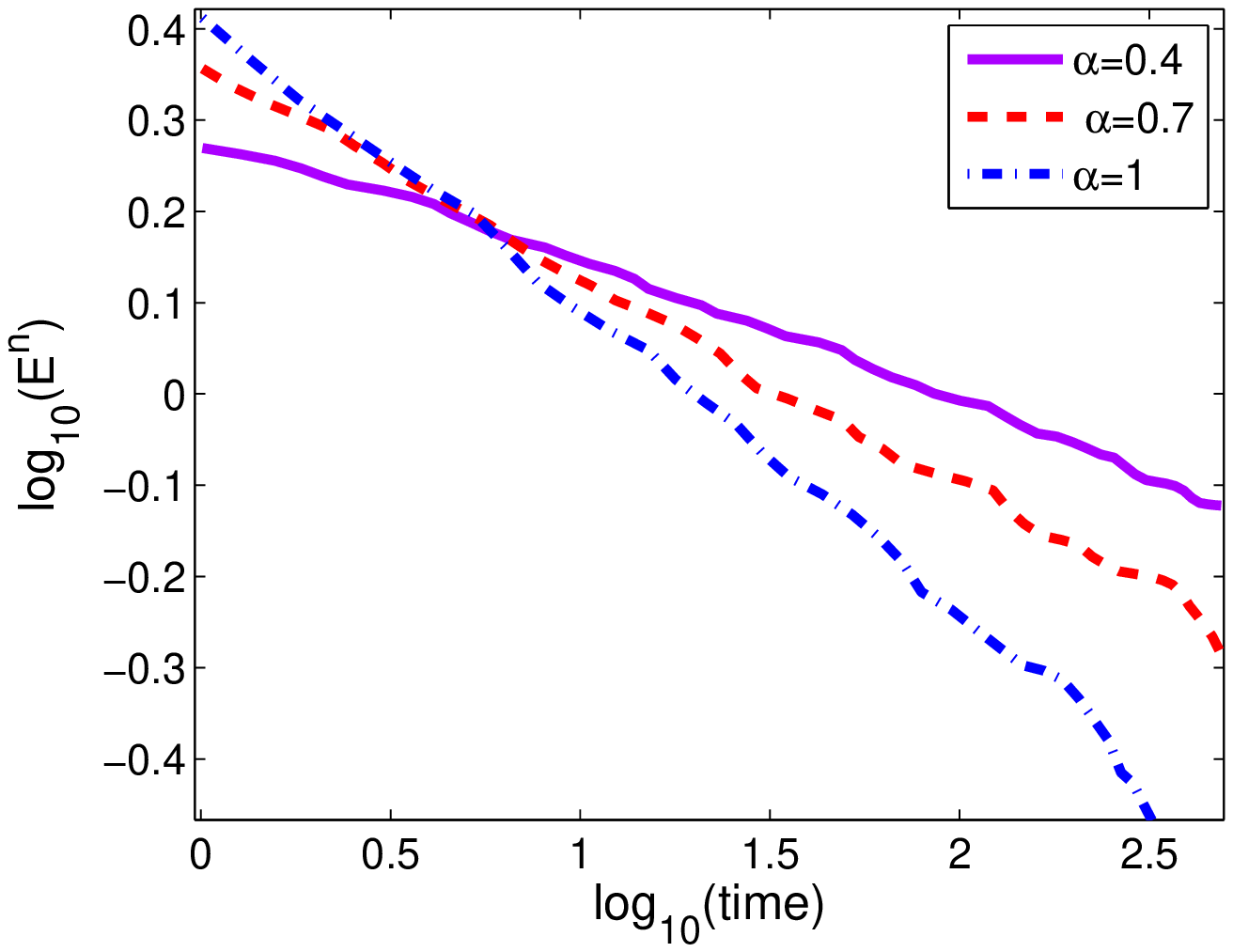}
\includegraphics[width=3.0in,height=2.0in]{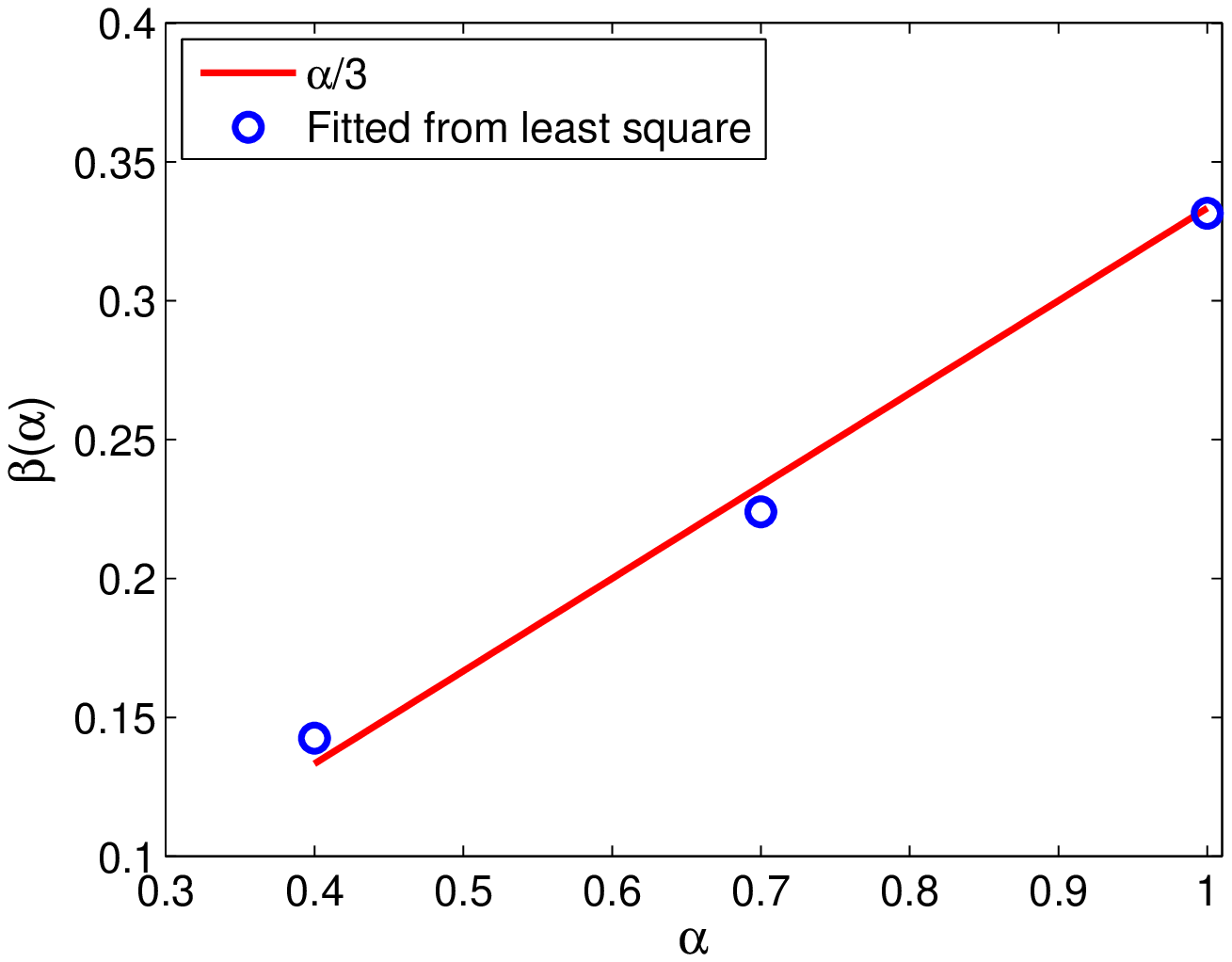}\\
\caption{The evolution of energy (left) and the least square fitted energy dissipation law scaling $\beta(\alpha)$ (right) of the Slope-Model for fractional order $\alpha=0.4,0.7$ and $1$, respectively.}
\label{Slope-Dynamic-Energy-Power}
\end{figure}

\begin{figure}[htb!]
\centering
\includegraphics[width=3.0in,height=2.0in]{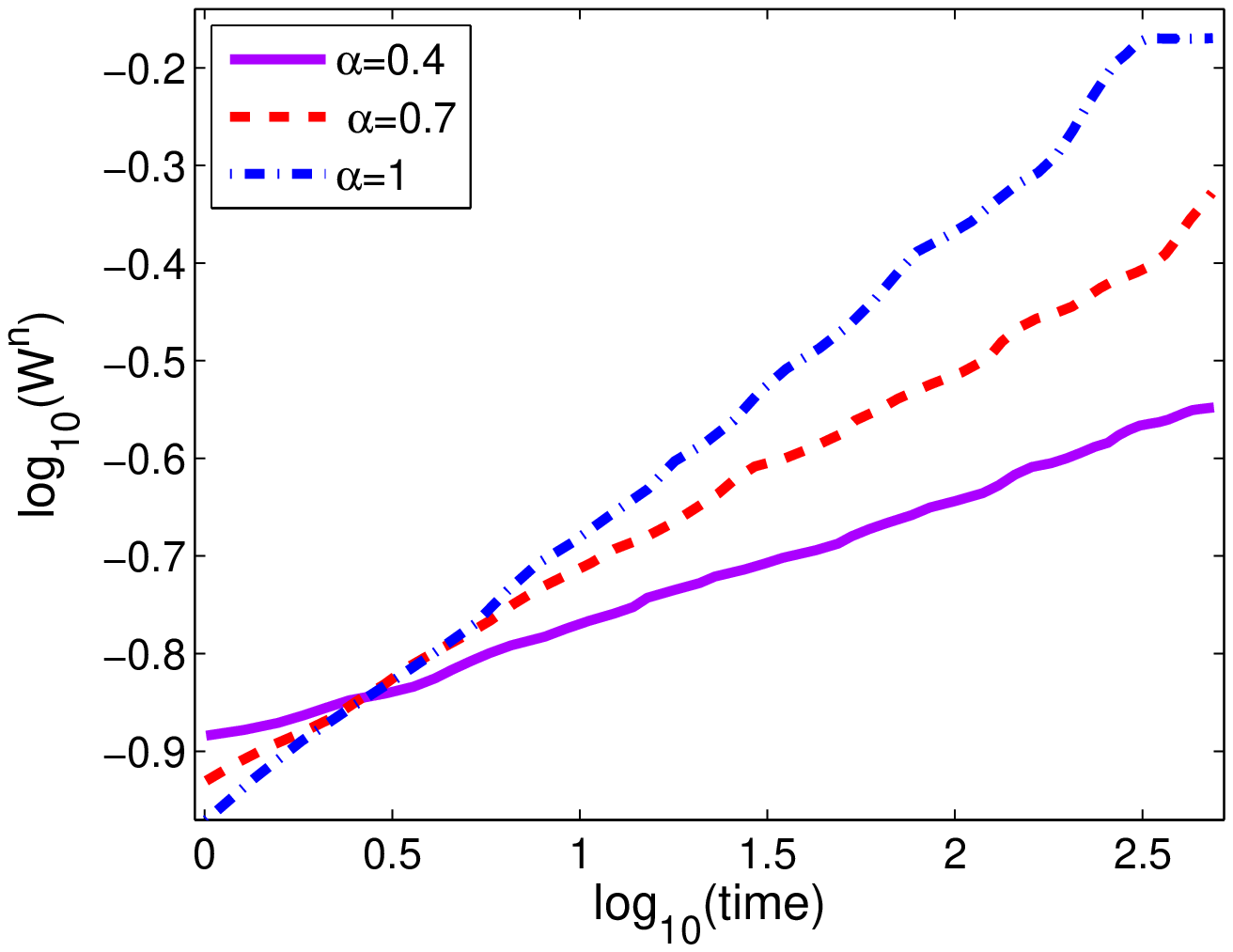}
\includegraphics[width=3.0in,height=2.0in]{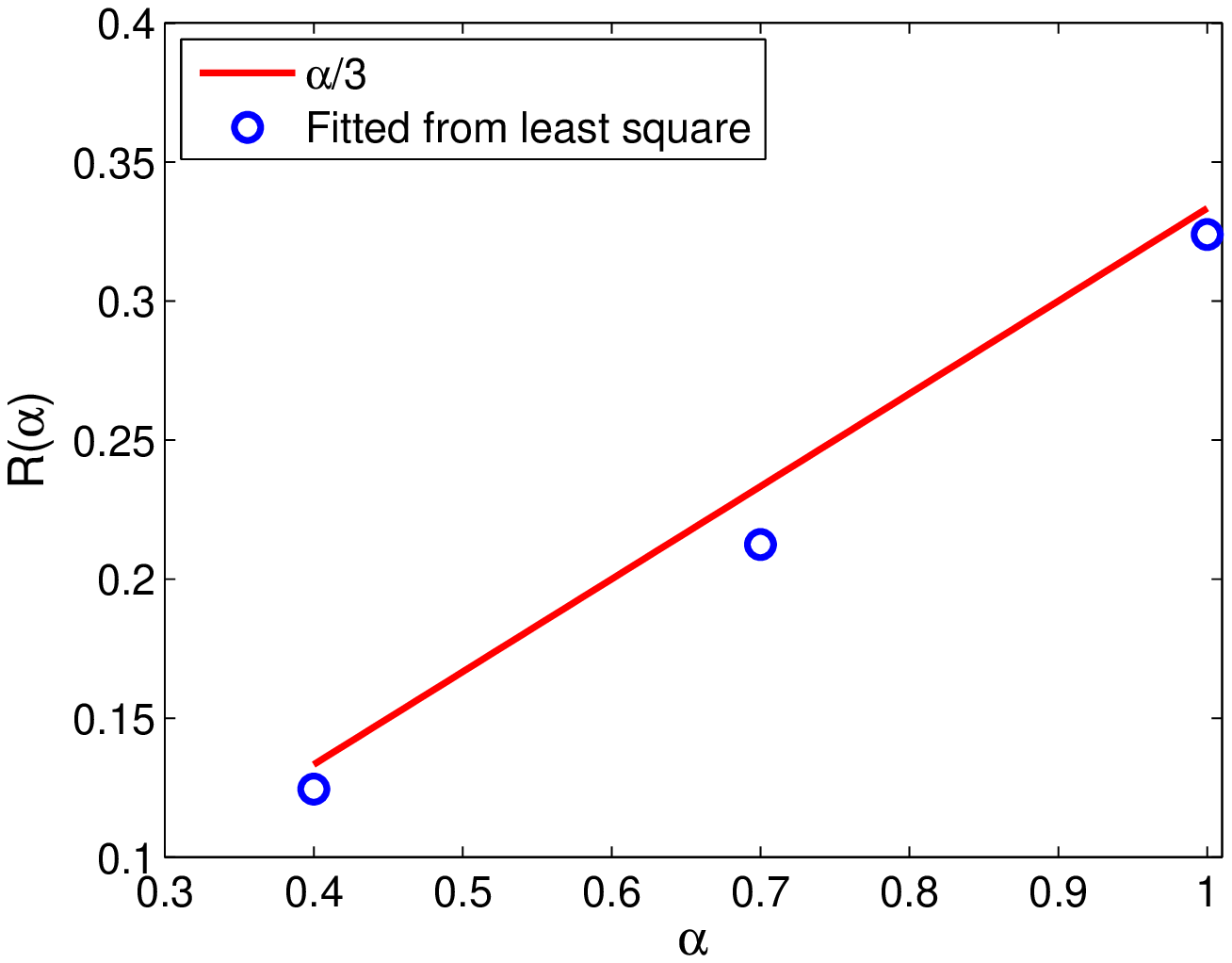}\\
\caption{The evolution of roughness (left) and the least square fitted roughness growth law scaling $R(\alpha)$ (right) of the Slope-Model for fractional order $\alpha=0.4,0.7$ and $1$, respectively.}
\label{Slope-Dynamic-Roughness-Power}
\end{figure}

Figure \ref{Slope-Model-Adaptive-Test} shows the plots of the energies and roughnesses for the Slope-Model and No-Slope-Model using uniform,
graded and adaptive time steps, respectively.
In our computation, we see that using uniform time step may produce incorrect steady-state solution while the energies still decay monotonically
and the roughnesses show analogous plots compared with the reference solution.
However, a comparison between the numerical results of using graded and adaptive time steps, we see that the curves are indistinguishable.
For the adaptive strategy, the density of circle indicates the size of the time step, we then also observe that small time steps are used at the early stage of the computation because the quick transition of solution.
Subsequent, large time steps are allowed due to the solution changes slowly.
Again, small time steps are employed when capturing the steep structural transition from one stage to the next one.
As a result, the total numbers of adaptive time steps are 4443 and 3586 for the Slope-Model and No-Slope-Model, respectively, while it takes 30000 constant time steps.
The above observations show that the effect of the adaptive time approach on efficiency is significant dramatic.

Below we make a comparison among different fractional orders $\alpha$ to exploit how fractional index affects the evolution dynamics.
Always, the third time mesh strategy is employed to solve the problem \eqref{Problem-2} with initial condition \eqref{Given-Initial-Condition} in what follows.
Figures \ref{Slope-Model-Energy-Roughness} and \ref{No-Slope-Model-Energy-Roughness} are the plots of the energies and roughnesses of the Slope-Model and No-Slope-Model for three different values of $\alpha$.
We observe that in all cases the original energy decays rapidly and the smaller fractional order $\alpha$ the faster the energy dissipates, later it decays slower as smaller fractional index $\alpha$, and they reach the analogous steady-state in the end.
The above observations may indicate that the time-fractional operator could affects the time scaling of the evolution dynamics, while the steady-state may not be affected.

\begin{example}\label{Coarsening-Dynamics}
We investigate here the coarsening dynamics using the CN-SAV schemes \eqref{Slope-CN-SAV-1}-\eqref{Slope-CN-SAV-2} and \eqref{No-CN-SAV-1}-\eqref{No-CN-SAV-2}.
The initial condition is a random state by assigning a random number varying from $-0.001$ to $0.001$ to each grid points with parameter value $\epsilon=0.03$.
\end{example}

To discover the scaling of effective free energy and roughness in the time-fractional MBE model during coarsening, define the absolute value of the Slope-Model for each linear line of energy and roughness, $\beta(\alpha),R(\alpha)$, as
\begin{align}
\log_{10}\bra{E(\alpha,t)}=\beta^{0}(\alpha)-\beta(\alpha)\log_{10}(t),\quad
\log_{10}\bra{W(\alpha,t)}=R^{0}(\alpha)+R(\alpha)\log_{10}(t),
\end{align}
and $E(\alpha,t)=\beta^{0}(\alpha)-\beta(\alpha)\log_{10}(t)$ for the No-Slope-Model, where $E(\alpha,t)$ and $W(\alpha,t)$ correspond to the energy $E(t)$ and roughness $W(t)$ with the fractional index $\alpha$, respectively.

\begin{figure}[htb!]
\centering
\includegraphics[width=3.0in,height=2.0in]{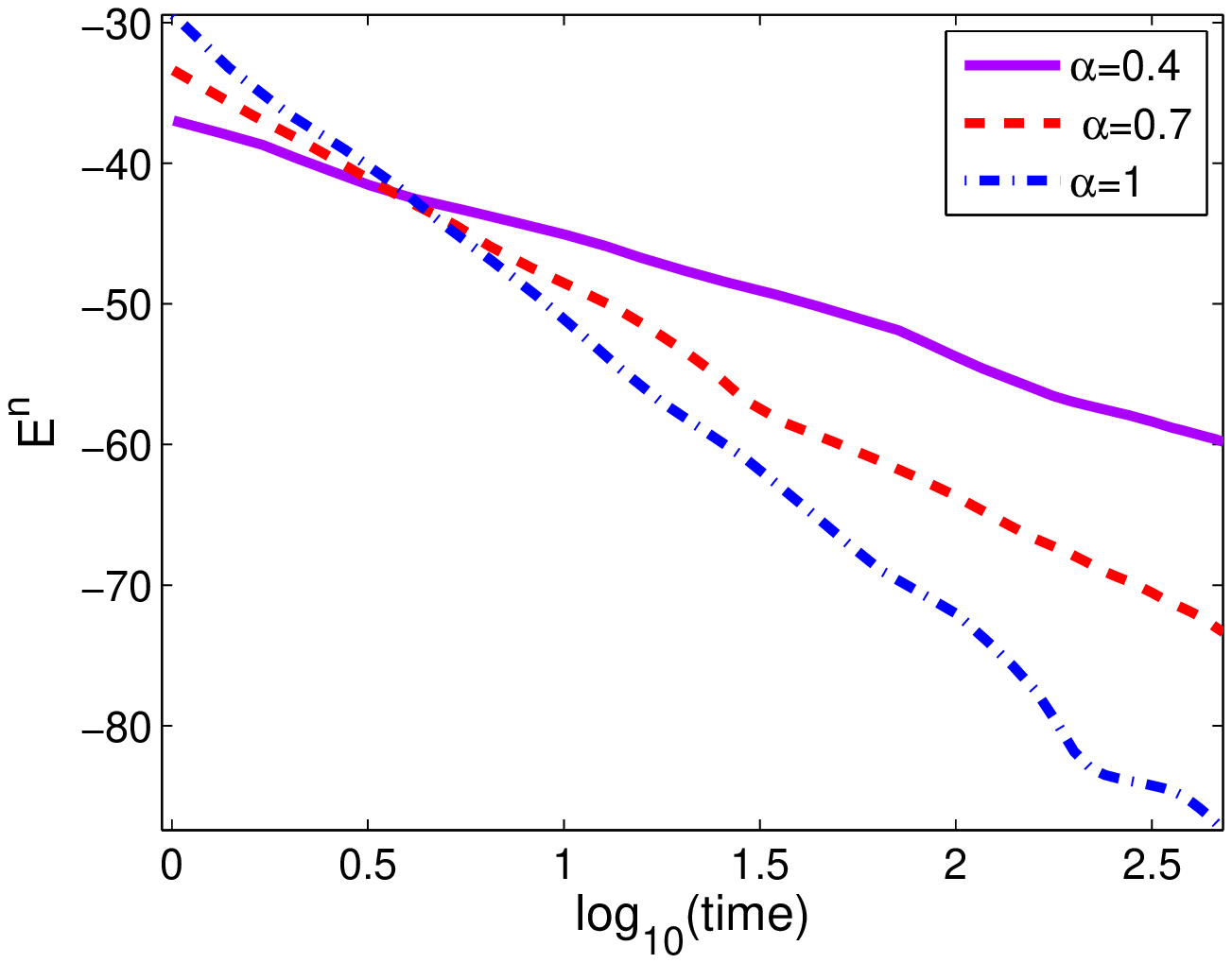}
\includegraphics[width=3.0in,height=2.0in]{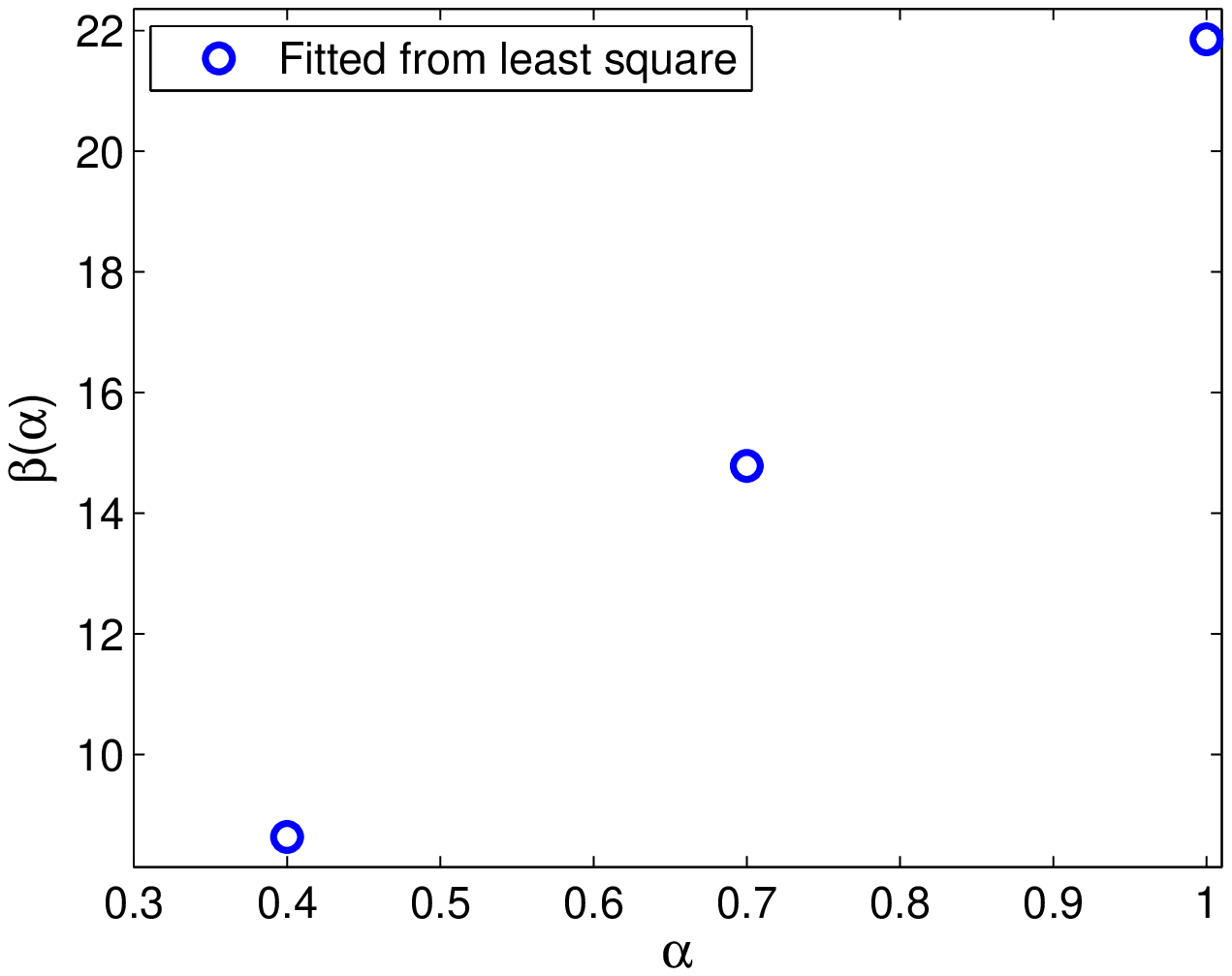}\\
\caption{The evolution of energy (left) and the least square fitted energy dissipation law scaling $\beta(\alpha)$ (right) of the No-Slope-Model for fractional order $\alpha=0.4,0.7$ and $1$, respectively.}
\label{No-Slope-Dynamic-Energy-Power}
\end{figure}

\begin{figure}[htb!]
\centering
\includegraphics[width=3.0in,height=2.0in]{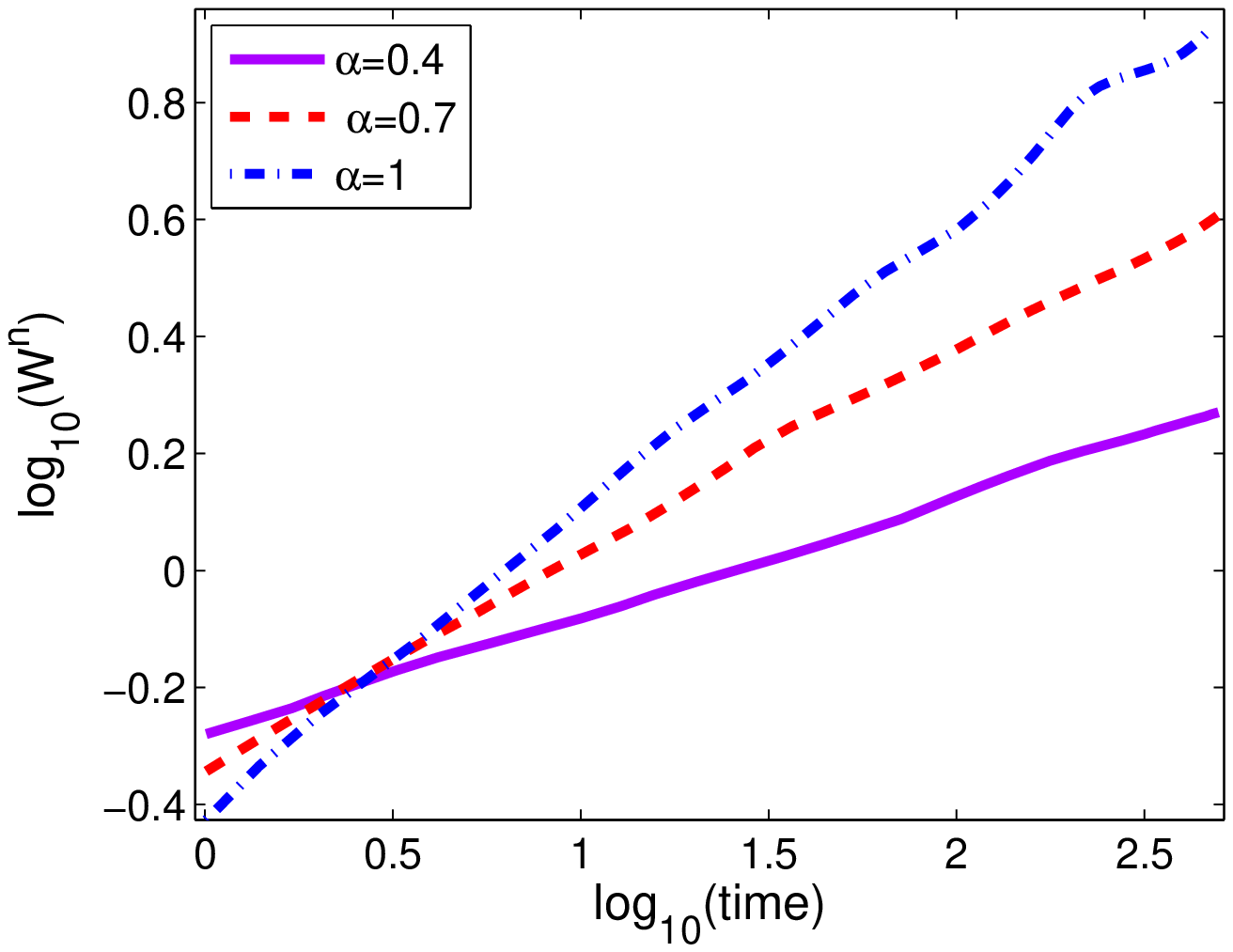}
\includegraphics[width=3.0in,height=2.0in]{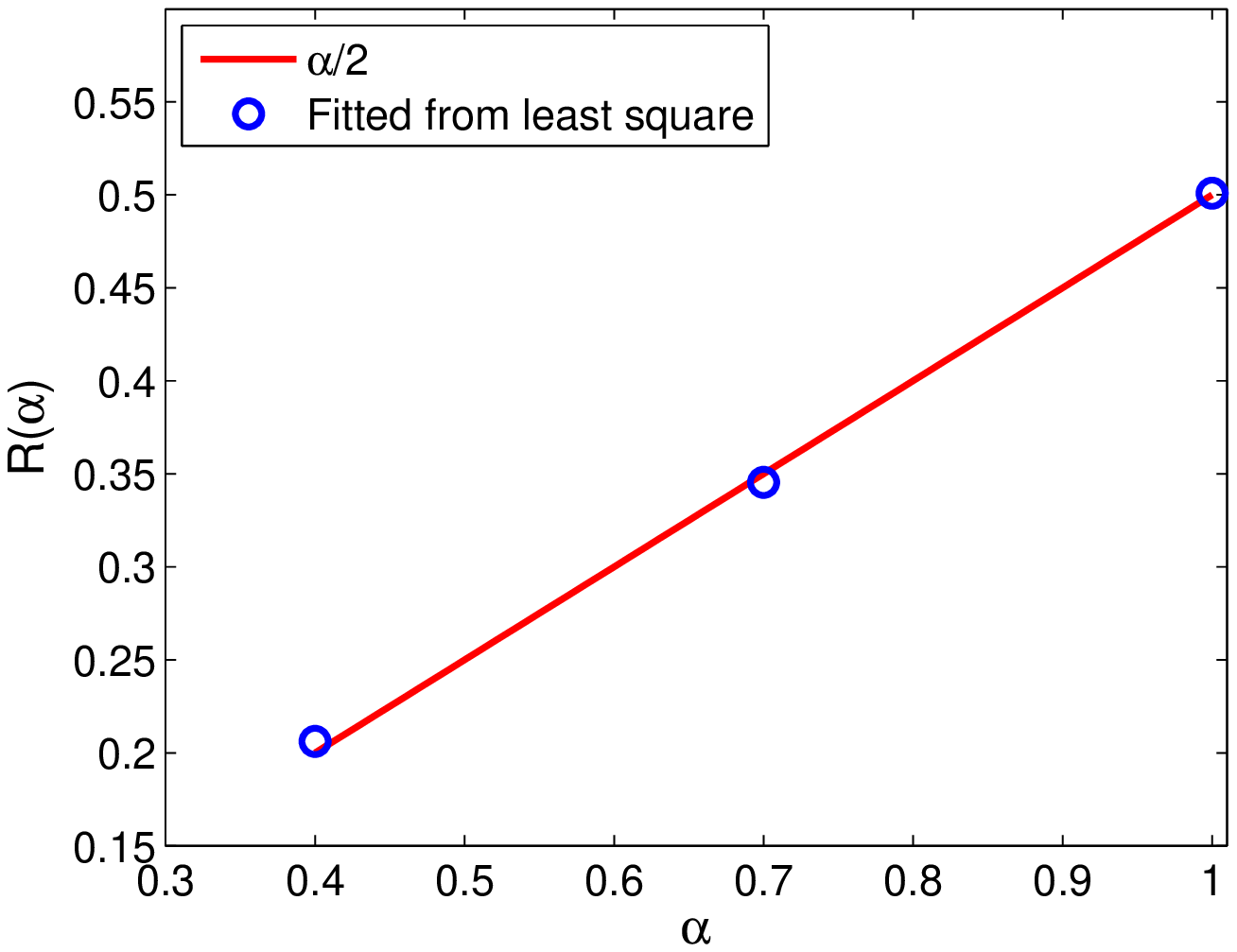}\\
\caption{The evolution of roughness (left) and the least square fitted roughness growth law scaling $R(\alpha)$ (right) of the No-Slope-Model for fractional order $\alpha=0.4,0.7$ and $1$, respectively.}
\label{No-Slope-Dynamic-Roughness-Power}
\end{figure}

Figures \ref{Slope-Dynamic-Energy-Power} and \ref{Slope-Dynamic-Roughness-Power} show the time evolutions of energy and roughness with different fractional orders $\alpha$ for the Slope-Model during $t\in[1,500]$
where the adaptive time-stepping parameters $\tau_{\min}=\tau_{N_{0}}=1.25\times10^{-4}$
and $\tau_{\max}=10^{-1}$.
We observe that
the energy dissipation approximately as $O(t^{-\frac{\alpha}{3}})$
and
the growth rate of the roughness is approximately as $O(t^{\frac{\alpha}{3}})$
which are consistent with consistent with $O(t^{-\frac{1}{3}})$ and $O(t^{\frac{1}{3}})$, respectively, as recorded in
\cite{Yang2017numerical,Cheng2019Highly}
for the integer Slope-Model.
The time evolutions of energy and roughness for the No-Slope-Model using adaptive parameters $\tau_{\min}=\tau_{N_{0}}=3.32\times10^{-5}$
and $\tau_{\max}=10^{-1}$ that
are depicted in Figures \ref{No-Slope-Dynamic-Energy-Power} and \ref{No-Slope-Dynamic-Roughness-Power}, respectively.
We see that the energy dissipation approximately as $O(-\log_{10}(t))$ and the growth rate of the roughness is approximately as $O(t^{\frac{\alpha}{2}})$ that are consistent with consistent with $O(-\log_{10}(t))$ and $O(t^{\frac{1}{2}})$ as reported in
\cite{Yang2017numerical,Cheng2019Highly}
for the integer No-Slope-Model, respectively.

Furthermore, the drawings of the Figure \ref{Slope-Dynamics-alpha-04-07-10}
displays the numerical solutions of the height function $\phi$ and its Laplacian $\Delta\phi$ for the Slope-Model with different fractional orders $\alpha$.
Based on Figure \ref{Slope-Dynamics-alpha-04-07-10} and additional results not shown here for brevity,
we observe that the edges of the pyramids generate a random distributed network over the domain and the pyramids become large when time increases.
Additionally, coarsening dynamics, at beginning, appear to be faster as smaller $\alpha$ while it would be much slower as time evolves.
Also, the observed phenomena are in good agreement with the published results \cite{Zhao2019On}.


\section{Concluding remarks}
In simulating the time-fractional phase field equations including the Molecular Beam Epitaxial model considered in this paper,
the initial singularity should be treated properly because it always destroys the time accuracy of numerical algorithms especially
near the initial time.

However, it seems challenging to build time-stepping approaches maintaining the discrete energy dissipation law based on the conventional L1 formula,
especially on general nonuniform time meshes.
Nonetheless, the energy stable schemes permitting adaptive time-stepping strategies are very attractive because
they could be applicable for time-fractional phase-field models
and for long-time simulations approaching the steady state.
As an interesting remedy, the novel L1$^{+}$ formula, is proposed to approximate the fractional Caputo's derivative.
As a consequence, coupled with  the scalar auxiliary variable, we suggest two linearized second-order energy stable CN-SAV schemes for the MBE model with and without slope selection, respectively,  by virtue of the naturally positive semi-definite property of a discrete quadratic form.
Furthermore, for the long-time simulations approaching the steady state,  the fast L1$^{+}$ version  incorporated with adaptive time-stepping strategy is developed for time-fractional phase field equations.
Ample numerical examples are presented to validate the effectiveness of CN-SAV schemes.

The L1$^{+}$ formula would be superior to some widespread approximations, such as the
L1, Alikhanov and BDF2-like formulas, because it is second-order accurate for both smooth and non-smooth solutions, and the convergence order is independent of the fractional order $\alpha\in(0,1)$.
Thus it becomes critical to establish the rigorous theory on consistency, stability
and convergence of nonuniform L1$^{+}$ formula.
These issues will be addressed in a forthcoming report.


\begin{thebibliography}{10}

\bibitem{Clarke1987Origin}
S.~Clarke and D.~Vvedensky.
\newblock Origin of reflection high-energy electron-diffraction intensity
  oscillations during molecular-beam epitaxy: a computational modeling
  approach.
\newblock {\em Phys. Rev. Lett.}, 58:2235--2238, 1987.

\bibitem{Villain1991Continuum}
J.~Villain.
\newblock Continuum models of crystal growth from atomic beams with and without
  desorption.
\newblock {\em Journal De Physique I}, 1:19--42, 1991.

\bibitem{Gyure1998Level}
M.~Gyure, C.~Ratsch, B.~Merriman, R.~Caflisch, S.~Osher, J~Zinck, and
  D.~Vvedensky.
\newblock Level-set methods for the simulation of epitaxial phenomena.
\newblock {\em Phys. Rev. E}, 58:6927--6930, 1998.

\bibitem{Shen2012Second}
J.~Shen, C.~Wang, X.~Wang, and S.~Wise.
\newblock Second-order convex splitting schemes for gradient flows with
  ehrlich-schwoebel type energy: application to thin film epitaxy.
\newblock {\em SIAM J. Numer. Anal.}, 50:105--125, 2012.

\bibitem{Xu2006Stability}
C.~Xu and T.~Tang.
\newblock Stability analysis of large time-stepping methods for epitaxial
  growth models.
\newblock {\em SIAM J. Numer. Anal.}, 44:1759--1779, 2006.

\bibitem{Yang2017numerical}
X.~Yang, J.~Zhao, and Q.~Wang.
\newblock Numerical approximations for the molecular beam epitaxial growth
  model based on the invariant energy quadratization method.
\newblock {\em J. Comput. Phys.}, 333:104--127, 2017.

\bibitem{Cheng2019Highly}
Q.~Cheng, J.~Shen, and X.~Yang.
\newblock Highly efficient and accurate numerical schemes for the epitaxial
  thin film growth models by using the {SAV} approach.
\newblock {\em J. Sci. Comput.}, 78:1467--1487, 2019.

\bibitem{Gong2019Energy}
Y.~Gong and J.~Zhao.
\newblock Energy-stable {Runge-Kutta} schemes for gradient flow models uing the
  energy quadratization approach.
\newblock {\em Appl. Math. Lett.}, 94:224--231, 2019.

\bibitem{Zheng2017A}
Z.~Li, H.~Wang, and D.~Yang.
\newblock A space-time fractional phase-field model with tunable sharpness and
  decay behavior and its efficient numerical simulation.
\newblock {\em J. Comput. Phys.}, 347:20--38, 2017.

\bibitem{Hou2017Numerical}
T.~Hou, T.~Tang, and J.~Yang.
\newblock Numerical analysis of fully discretized {C}rank-{N}icolson scheme for
  fractional-in-space {Allen-Cahn} equations.
\newblock {\em J. Sci. Comput.}, 72:1--18, 2017.

\bibitem{Liu2018Time}
H.~Liu, A.~Cheng, H.~Wang, and J.~Zhao.
\newblock Time-fractional {A}llen-{C}ahn and {C}ahn-{H}illiard phase-field
  models and their numerical investigation.
\newblock {\em Comp. Math. Appl.}, 76:1876--1892, 2018.

\bibitem{Zhao2019On}
J.~Zhao, L.~Chen, and H.~Wang.
\newblock On power law scaling dynamics for time-fractional phase field models
  during coarsening.
\newblock {\em Commu. Non. Sci. Numer. Simul.}, 70:257--270, 2019.

\bibitem{Tang2018On}
T.~Tang, H.~Yu, and T.~Zhou.
\newblock On energy dissipation theory and numerical stability for
  time-fractional phase field equations.
\newblock {\em arXiv:1808.01471v1}, 2018.

\bibitem{Jin2016An}
B.~Jin, R.~Lazarov, and Z.~Zhou.
\newblock An analysis of the {L1} scheme for the subdiffusion equation with
  nonsmooth data.
\newblock {\em IMA J. Numer. Anal.}, 36:197--221, 2016.

\bibitem{JinLiZhou:2017}
B.~Jin, B.~Li, and Z.~Zhou.
\newblock Numerical analysis of nonlinear subdiffusion equations.
\newblock {\em SIAM J. Numer. Anal.}, 56:1--23, 2018.

\bibitem{Liao2018Sharp}
H.-L. Liao, D.~Li, and J.~Zhang.
\newblock Sharp error estimate of nonuniform {L}1 formula for time-fractional
  reaction-subdiffusion equations.
\newblock {\em SIAM J. Numer. Anal.}, 56:1112--1133, 2016.

\bibitem{Liao2018Unconditional}
H.-L. Liao, Y.~Yan, and J.~Zhang.
\newblock Unconditional convergence of a fast two-level linearized algorithm
  for semilinear subdiffusion equations.
\newblock {\em J. Sci. Comput.}, 80(1) (2019), 1-25.

\bibitem{Liao2018second}
H.-L. Liao, W.~Mclean, and J.Zhang.
\newblock A second-order scheme with nonuniform time steps for a linear
  reaction-sudiffusion problem.
\newblock {\em arXiv:1803.09873v4}, 2019.

\bibitem{Mclean1996Discretization}
W.~Mclean, V.~Thom\'ee, and L.~Wahlbin.
\newblock Discretization with variable time steps of an evolution equation with
  a positive-type memory term.
\newblock {\em J. Comp. Appl. Math.}, 69:49--69, 1996.

\bibitem{William2007A}
W.~McLean and K.~Mustapha.
\newblock A second-order accurate numerical method for a fractional wave
  equation.
\newblock {\em Numer. Math.}, 105:481--510, 2007.

\bibitem{Alikhanov2015A}
A.~Alikhanov.
\newblock A new difference scheme for the time fractional diffusion equation.
\newblock {\em J. Comput. Phys.}, 280:424--438, 2015.

\bibitem{Liao2016A}
H.-L. Liao, Y.~Zhao, and X.~Teng.
\newblock A weighted {ADI} scheme for subdiffusion equations.
\newblock {\em J. Sci. Comput.}, 69:1144--1164, 2016.

\bibitem{Gao2014A}
G.~Gao, Z.~Sun, and H.~Zhang.
\newblock A new fractional numerical differentiation formula to approximate the
  {C}aputo fractional derivative and its applications.
\newblock {\em J. Comput. Phys.}, 259:33--50, 2014.

\bibitem{LvXu2016}
C.~Lv and C.~Xu.
\newblock Error analysis of a high order method for time-fractional diffusion
  equations.
\newblock {\em SIAM J. Sci. Comput.}, 38:A2699--A2724, 2016.

\bibitem{Liao2016Stability}
H.-L. Liao, P.~Lyu, S.~Vong, and Y.~Zhao.
\newblock Stability of fully discrete schemes with interpolation-type
  fractional formulas for distributed-order subdiffusion equations.
\newblock {\em Numer. Algo.}, 75:845--878, 2017.

\bibitem{Liao2018discrete}
H.-L. Liao, W.~Mclean, and J.~Zhang.
\newblock A discrete {G}r\"{o}nwall inequality with applications to numerical
  schemes for subdiffusion problems.
\newblock {\em SIAM J. Numer. Anal.}, 57:218--237, 2019.

\bibitem{Jiang2017Fast}
S.~Jiang, J.~Zhang, Z.~Qian, and Z.~Zhang.
\newblock Fast evaluation of the {C}aputo fractional derivative and its
  applications to fractional diffusion equations.
\newblock {\em Commu. Comput. Phys.}, 21:650--678, 2017.

\bibitem{Gomez2011Provably}
H.~Gomez and T.~J. Hughes.
\newblock Provably unconditionally stable, second-order time-accurate, mixed
  variational methods for phase-field models.
\newblock {\em J. Comput. Phys.}, 230:5310--5327, 2011.

\bibitem{Qiao2011An}
Z.~Qiao, Z.~Zheng, and T.~Tang.
\newblock An adaptive time-stepping strategy for the molecular beam epitaxy
  models.
\newblock {\em SIAM J. Sci. Comput.}, 22:1395--1414, 2011.

\end{thebibliography}
\end{document}